\documentclass[12pt]{article}
\usepackage{amsfonts,amssymb,version}

\def\F{{\mathcal F}}
\def\GG{{\mathbb G}}

\def\hn{\mathop{\rm HN}\nolimits}

\def\mm{\mathfrak m}
\def\K{{\mathcal K}}
\def\L{{\mathcal L}}
\def\LL{{\mathbb L}} 

\def\OO{{\mathcal O}}

\def\Q{\mathbb Q}

\def\Z{\mathbb Z}

\def\deg{\mathop{\rm deg}\nolimits} 
\def\depth{\mathop{\rm depth}\nolimits}
 
\def\hom{\mathop{\underline{Hom}}\nolimits}

\def\Proj{\mathop{\rm Proj}\nolimits}

\def\rank{\mathop{\rm rank}\nolimits}
\def\spec{\mathop{\rm Spec}\nolimits}

\def\Sym{\mathop{\rm Sym}\nolimits}
\def\sym{\mathop{\rm Sym}\nolimits}

\let\hra\hookrightarrow
\let\ov\overline
\let\un\underline
\let\wh\widehat

\newtheorem{theorem}{Theorem}[section]
\newtheorem{proposition}[theorem]{Proposition}
\newtheorem{lemma}[theorem]{Lemma}

\newtheorem{corollary}[theorem]{Corollary}

\def\stm{\refstepcounter{theorem}\paragraph{\thetheorem}}
\def\rem{\refstepcounter{theorem}\paragraph{Remark \thetheorem}}

\def\defin{\refstepcounter{theorem}\paragraph{Definition \thetheorem}}
\makeatletter \def\l@section{\@dottedtocline{1}{0em}{1.2em}} \makeatother

\begin{document}

\centerline{\Large\bf Schematic Harder-Narasimhan stratification for}

\centerline{\Large\bf families of principal bundles in higher dimensions}

\bigskip

\centerline{\bf Sudarshan Gurjar and Nitin Nitsure} 


\begin{abstract}
Let $G$ be a connected split reductive group over a
field $k$ of characteristic zero. 
Let $X\to S$ be a smooth projective morphism of $k$-schemes, 
with geometrically connected fibers. We formulate a 
natural definition of a relative canonical reduction, 
under which principal 
$G$-bundles of any given Harder-Narasimhan type $\tau$ on fibers of 
$X/S$ form an Artin algebraic stack $Bun_{X/S}^{\tau}(G)$ over $S$, 
and as $\tau$ varies, these 
stacks define a stratification of the stack $Bun_{X/S}(G)$ by 
locally closed substacks.
This result extends to principal bundles 
in higher dimensions the earlier such result
for principal bundles on families of curves.
The result is new even for vector bundles, that is, for $G = GL_{n,k}$.
\end{abstract}


\centerline{2010 Math. Subj. Class. : 14D20, 14D23, 14F10.}


\section{Introduction}
Let $G$ be a connected split reductive algebraic group over a 
field $k$ of characteristic zero. Let there be chosen 
a split maximal torus and a Borel containing it, 
and let $\ov{C}$ 
denote the corresponding closed positive Weyl chamber.
If $(X,\OO_X(1))$ is a projective variety over an extension field $K/k$ 
together with a very ample line bundle, 
and if $E$ is a principal $G$-bundle on $X$,
then recall that the Harder-Narasimhan type of $E$ is an element 
$\hn(E) \in \ov{C}$. It is defined as the type
of the canonical reduction of $E$, which is a particular
rational reduction of the structure group of $E$, 
from $G$ to a standard parabolic $P$. Here, a 
rational reduction is a reduction defined  
on a big open subscheme of $X$, that is, 
an open subscheme whose complement is of codimension $\ge 2$.

We now move to the relative set-up. Let $X\to S$ be a 
smooth projective morphism with geometrically connected fibers, 
where $S$ is a noetherian scheme over $k$, with  
a given relatively very ample line bundle.
A family of principal $G$-bundles on $X/S$ means 
a principal $G$-bundle $E$ on $X$, so that
each restriction $E_s = E|X_s$ is a principal $G$-bundle on the 
smooth projective variety
$X_s$ which is the fiber of $X$ over $s\in S$. The individual
Harder-Narasimhan types $\hn(E_s)\in \ov{C}$ together define a 
function $S \to \ov{C}$. We prove that (see Proposition 
\ref{semicontinuity and uniqueness} assertion (1))
this function is upper semi-continuous
w.r.t. the usual partial ordering on $\ov{C}$.  This implies that for 
any $\tau \in \ov{C}$, the subset $|S|^{\le \tau}(E)$ 
which consists
of all $s\in S$ such that $\hn(E|X_s)\le \tau$ is open in $S$, 
the subset $|S|^{\tau}(E)$ which consists
of all $s\in S$ such that $\hn(E|X_s)= \tau$ is closed in $|S|^{\le \tau}(E)$, 
and the closure of $|S|^{\tau}(E)$ in $|S|$ is contained in 
$\bigcup_{\gamma \ge \tau}\, |S|^{\gamma}(E)$. 
Hence for each $\tau$ the Artin stack $Bun_{X/S}(G)$ of principal 
$G$-bundles over
$X/S$ has an open substack $Bun_{X/S}^{\le \tau}(G)$ which to any $S$-scheme $T$ 
associates the groupoid whose objects are all principal $G$-bundles
$E$ on $X_T$ such that $\hn(E_t) \le \tau$ at all $t\in T$.

If we are interested not just in the discrete 
invariant $\hn(E_s)$ but in the behaviour of the canonical reductions
of $E_s$ in a family, then we need to have a good, workable
definition of a relative rational reduction of structure group
of such an $E$, from $G$ to a standard parabolic $P$. On the one hand,
in the special case where $S = \spec K$ for an extension field $K/k$,
the definition of a relative rational reduction to $P$ 
should amount just to a reduction to $P$ over a big open subscheme 
$U\subset X$. On the other hand, the definition should have a modicum of 
$\OO$-coherence and 
flatness built into it which would allow us to deploy the theory
of flatness and base change for coherent sheaves and their cohomologies.


\pagestyle{myheadings}
\markright{Gurjar and Nitsure: Principal bundles
in higher dimensions.}


We now describe our candidate for such a definition. To begin with,
to each standard parabolic $P$, we associate the irreducible 
linear $G$-representation $V_P = \Gamma(G/P,\omega^{-1}_{G/P})^{\vee}$. 
As the anti-canonical
line bundle $\omega^{-1}_{G/P}$ is very ample, it
defines a $G$-equivariant embedding  $G/P \hra {\bf P}(V_P)$, which denotes the 
projective space of lines in $V_P$. 

Let $X\to S$ be a 
smooth projective morphism with geometrically connected fibers, 
where $S$ is a noetherian scheme over $k$, with  
a given relatively very ample line bundle $\OO_{X/S}(1)$ on $X$.
Let $E$ be a principal $G$-bundle on $X$.  
Let $E(V_P)$ denote the vector bundle associated to $E$ by 
the $G$-representation $V_P$.
We define a {\bf relative rational reduction} of structure group
of $E$ from $G$ to $P$ to be a pair $(L,f)$, 
where $L$ is a line bundle on $X$ 
and $f: L\to E(V_P)$ is an injective $\OO_X$-linear 
homomorphism of sheaves,
such that \\
(i) the open subscheme $U = \{ x\in X\,|\,\rank(f_x) =1\} \subset X$ 
is relatively big over $S$, that is, for each $s\in S$
the fiber $U_s$ has complementary codimension 
$\ge 2$ in the fiber $X_s$, and \\
(ii) the section $U\to {\bf P}(E(V_P))$  
defined by $f$ factors via the natural closed embedding 
$E/P \hra {\bf P}(E(V_P))$. 

Note that a section $\sigma: 
U\to {\bf P}(E(V_P))$ is the same as a line subbundle 
$f': L'\hra E(V_P)|U$, which is the pullback by $\sigma$ 
of the tautological line subbundle $\OO(-1) \hra E(V_P)_{{\bf P}(E(V_P))}$.
By Proposition \ref{line bundle must admit a prolongation}, giving 
the extra data $(L,f)$ simply amounts to 
imposing the requirement that $L'$
should admit a prolongation to a line bundle $L$ on $X$. If it exists,
such a prolongation $L$ will be unique, and in that case there will 
exist a unique map $f: L \to E(V_P)$ which prolongs $f'$.

In the special case where $S = \spec K$ for a field $K$,
the above definition is equivalent to the usual definition 
(this is the Proposition \ref{rational parabolic reductions}), 
as one would require of any such generalization. 
Finally, we define a {\bf relative canonical reduction} for $E$ on $X$
over $S$ to be a relative rational reduction $(L,f)$ of the structure group
of $E$ from $G$ to a standard parabolic $P$, such that for each $s\in S$,
the restriction $(L|X_s, f|X_s)$ is a canonical reduction of $E_s = E|X_s$.

For any type $\tau \in \ov{C}$, the above definition allows us to define
an $S$-groupoid $Bun_{X/S}^{\tau}(G)$, which attaches to any $S$-scheme $T$ 
the category $Bun_{X/S}^{\tau}(G)(T)$ whose objects are $(E, L, f)$ where
$E$ is a principal $G$-bundle on $X_T$ and $(L,f)$ is a 
relative canonical reduction
for $E$ of constant type $\tau$. We prove that (see Proposition 
\ref{semicontinuity and uniqueness} assertion (2)) 
if a relative canonical reduction exists, then 
it is unique up to a unique isomorphism, that is, 
given
any two such reductions $(L,f)$ and $(L',f')$, there exists 
a unique isomorphism $\phi: L\to L'$ such that $f = f'\circ \phi$.  
The following is our main result.
 
\begin{theorem}\label{stacky main theorem} 
For any $\tau \in \ov{C}$, the 
forgetful $1$-morphism $Bun_{X/S}^{\tau}(G) \to Bun_{X/S}(G)$ is 
a locally closed embedding of Artin stacks. 
As $\tau$ varies over $\ov{C}$, this defines a 
stratification of $Bun_{X/S}(G)$ by the locally closed substacks 
$Bun_{X/S}^{\tau}(G)$
with respect to the  standard partial order on $\ov{C}$, that is,
each $Bun_{X/S}^{\tau}(G)$ is a closed substack of the open substack
$Bun_{X/S}^{\le \tau}(G)$ of $Bun_{X/S}(G)$. 
\end{theorem}

The reader is cautioned that while 
the closure of $Bun_{X/S}^{\tau}(G)$ in $ Bun_{X/S}(G)$ 
will be contained in $\cup_{\gamma \ge \tau}\,
Bun_{X/S}^{\gamma}(G)$, the inclusion will in general be proper.

The above result has the following formulation in elementary terms.


\begin{theorem}\label{main theorem} 
Let $E$ be a principal $G$-bundle
on a smooth projective family of varieties $X/S$. 
Then for each $\tau \in \ov{C}$, there exists a locally closed subscheme 
$S^{\tau}(E)\subset S$ with the following universal property:
A morphism $T\to S$ factors via the inclusion $S^{\tau}(E) \hra S$
if and only if the pullback $E_T$ of $E$ to $X_T$ admits a relative 
canonical reduction
of constant type $\tau$ over $T$. As $\tau$ varies over $\ov{C}$, 
this defines a 
stratification of $S$ by the locally closed subschemes $S^{\tau}(E)$
with respect to the  standard partial order on $\ov{C}$.
\end{theorem} 


Note that in particular, the Theorem \ref{main theorem} asserts that 
there exists a relative canonical reduction in our sense 
over each $S^{\tau}(E)$. This is stronger than just the conclusion
that there exists a relatively big open subscheme in 
$X\times_S S^{\tau}(E)$ over which we have a parabolic reduction
which restricts to the canonical reduction on each fiber $X_s$ 
for $s\in S^{\tau}(E)$.

When $S$ is reduced and $\hn(E_s)$ is constant over $S$, 
the above theorem gives the existence of a relative canonical reduction
over $S$, in the strong sense of our above new definition of such
a relative reduction. By forgetting the extra data involved,  
this has the following immediate consequence, which 
-- as the reader may note -- 
makes no reference at all to our new definition 
of a relative canonical reduction.

\begin{corollary}\label{reduced base}
Let $E$ be a principal $G$-bundle
on a smooth projective family of varieties $(X, \OO_{X/S}(1))$  
on a noetherian $k$-scheme $S$. 
If $S$ is reduced and if the Harder-Narasimhan type
$\hn(E_s)$ is constant over $S$, then there exists 
a relative canonical reduction $(L,f)$ of $E$ over $S$.
In particular, there exists  
a parabolic reduction $\sigma: U\to E/P$ of $E$ that is defined 
on a relatively big open subscheme $U$ of $X$ over $S$, 
which restricts to the canonical reduction on each fiber
of $X\to S$.
\end{corollary}

Now some history. The concept of HN-type 
was introduced for vector bundles on curves by Harder and Narasimhan
in [H-N], where they make use of the  
corresponding set-theoretic stratification. 
The existence of a set-theoretic HN-stratification for families
of vector bundles (in the sense of $\mu$-semi-stability) 
over higher dimensional projective varieties was proved 
by Shatz in [Sh]. Simpson extended these results to pure $\OO$-coherent
sheaves (in the sense of Gieseker semistability), and 
also to pure $\OO$-coherent sheaves 
of $\Lambda$-modules in the sense of Deligne, 
in [Si]. The concept of
canonical reduction for principal bundles was introduced in the context
of curves, and its existence and uniqueness proved due to the efforts of 
Ramanathan [Ram], Atiyah and Bott [A-B], and Behrend [Be 1], [Be 2]. 

The set-theoretic HN-stratification for pure $\OO$-coherent sheaves
on projective schemes,
in the sense of Gieseker semistability, was elevated in [Ni 3] 
to a scheme-theoretic
stratification, where the schematic strata have the appropriate universal 
property, similar to that in Theorem \ref{main theorem} above.   
An analogous result for pure $\OO$-coherent sheaves 
of $\Lambda$-modules on projective schemes was proved in [Gu-Ni 1].

The paper [Gu-Ni 1] also considers families of principal bundles over 
families of curves in characteristic zero, with reductive structure group, 
and proves the analog of Theorem \ref{main theorem} above over curves. 
We subsequently learned that this
had in fact been earlier proved in the PhD thesis of Behrend [Be 1]. 
In the present paper, we prove the result in all dimensions. 
Our proof is by induction on the relative dimension of $X/S$,
which we can now assume is $\ge 2$. 

In the case of curves, a deformation theoretic argument shows that 
the stacks $Bun_{X/S}^{\tau}(G)$ are all reduced (which holds 
also in characteristic $p$ if the Behrend conjecture holds for $G$ over $k$). 
In higher dimensions, we do not know even in characteristic zero 
whether $Bun_{X/S}^{\tau}(G)$ are always reduced. 

\rem Though we have written this paper for characteristic zero, everything
goes through unchanged when $k$ and $G$ are such that the Behrend
conjecture holds for $G$ and for all Levi quotients 
$P/R_u(P)$ of standard parabolics $P$ of $G$.

\rem In the special case $G = GL_{n,k}$, the main results of this paper 
(that is, Theorems \ref{stacky main theorem}, \ref{main theorem}, 
Corollary \ref{reduced base}, Proposition \ref{semicontinuity and uniqueness})
can be regarded as results about relative HN-filtrations and schematic 
HN-stratifications for families of vector
bundles in the sense of $\mu$-semistability, 
unlike the corresponding earlier results proved in  
[Ni 3] which were in the sense of Gieseker semistability.
Thus, these results are new even in the case of  
vector bundles, that is, for $G = GL_{n,k}$.

\medskip

This article is arranged as follows. The sections 2 to 5 are devoted to setting
up the basics of relative canonical reductions and their restrictions to
relative divisors $Y\subset X$ over $S$. The section 6 proves a result 
on embedding of 
relative Picard schemes $Pic_{X/S}\to Pic_{Y/S}$ which is needed for
lifting relative canonical reductions from a 
relative divisor $Y\subset X$ over $S$ to all of $X$.
The main results, including Propositions \ref{semicontinuity and uniqueness}, 
\ref{lifting the relative canonical reduction} and
Theorem \ref{main theorem}, are proved in section 7.

{\bf Acknowledgements.} We thank the 
Center for Quantum Geometry of Moduli Spaces for their hospitality, 
supported by a Center of Excellence grant from the 
Danish National Research Foundation (DNRF95) and 
by a Marie Curie International Research Staff Exchange Scheme
Fellowship within the 7th European Union Framework Programme
(FP7/2007-2013) under grant agreement n° 612534, project MODULI - Indo
European Collaboration on Moduli Spaces. Sudarshan Gurjar would also
like to thank the ICTP and the TIFR for support during part of the work. 
We thank Najmuddin Fakhruddin for 
producing the Example \ref{due to Najmuddin} of a line bundle
defined on a relatively big open subscheme which does not
admit a global prolongation.


\section{Linearized rational sections of projective \\bundles}

For any vector bundle 
(that is, a coherent locally free sheaf) $F$ on  a scheme $X$,
let ${\bf P}(F) = \Proj_{\OO_X}\,\Sym_{\OO_X}(F^{\vee})$ denote the 
projective bundle of lines in the fibers of $F$. 
Recall that such an $F$ is called a linearization of the projective
bundle ${\bf P}(F)\to X$. This motivates the word `linearized' in the 
following definition.

\defin\label{groupoid is equivalent to set} 
Given a vector bundle $F$ on a scheme $X$,
an {\bf $F$-linearized rational section} of the projective bundle
${\bf P}(F)\to X$ of lines in the fibers of $F$
will mean a rank $1$ locally free subsheaf $[L,f]$ of $F$, 
represented by a pair $(L,f)$ consisting of 
a rank $1$ locally free sheaf $L$ on $X$ together with an injective
$\OO_X$-linear homomorphism of sheaves $f: L\to F$, 
such that the open subscheme
$U = \{ x\in X \,|\, \rank(f_x) = 1\}$ is {\bf big in $X$}, that is, 
$X - U$ is of codimension $\ge 2$ in $X$ at all points
(here, $f_x = f|\spec \kappa(x)$ where $\kappa(x) = \OO_{X,x}/\mm_x$ 
denotes the residue field at $x$). 
We denote by $R(F/X)$
the set of all such $F$-linearized rational sections
$[L,f]$. 

Let $X\to S$ be a morphism of schemes, and let $F$ be a 
vector bundle on $X$. An {\bf $F$-linearized relative rational section
w.r.t. $S$} of the projective bundle ${\bf P}(F)\to X$  
will mean a rank $1$ locally free subsheaf $[L,f]$ of $F$,
represented by a pair $(L,f)$ consisting of 
a rank $1$ locally free sheaf $L$ together with 
an injective $\OO_X$-linear homomorphism of sheaves
$f: L\to F$, such that the open subscheme
$U = 
\{ x\in X \,|\, \rank(f_x) = 1\}$ is {\bf relatively big in $X$ over $S$}, 
that is, for each $s\in S$, the closed subset 
$X_s - U_s$ is of codimension $\ge 2$ in $X_s$ at all points of $X_s$. 
We will denote by $R(F/X/S)$
the set of all such $F$-linearized relative rational sections.
It is clear that $R(F/X/S) \subset R(F/X)$.
When the base $S$ is of the form $\spec K$ for a field $K$, 
then we get $R(F/X/\spec K) = R(F/X)$.

By definition, two pairs $(L_1,f_1)$ and $(L_2,f_2)$ represent the 
same element $[L,f]$ in $R(F/X)$ in the absolute case 
(or in $R(F/X/S)$ in the relative case)  
if and only if there exists an $\OO_X$-linear
isomorphism $\phi : L_1\to L_2$ such that $f_1 = f_2\circ \phi$.
Such a $\phi$, if it exists, is unique.

\begin{proposition}\label{line bundle must admit a prolongation}
Let $\pi: X\to S$ be a smooth morphism where
$S$ is noetherian, and let $F$ be a vector bundle on $X$.
Let $\OO_F(-1) \subset F_{{\bf P}(F)}$ denote the tautological line subbundle
of the pullback of $F$ to ${\bf P}(F)$.
The there is a natural bijection between the following two sets.

(i) The set $R(F/X/S)$ of all $F$-linearized rational sections of ${\bf P}(F)$.

(ii) The set $\Sigma(F/X/S)$ which consists of all pairs $(U,\sigma)$ where 
$U\subset X$ is an open subscheme which is relatively big over $S$ and
$\sigma$ is section of ${\bf P}(F)$ over $U$ 
which is maximal in the sense that $\sigma$ cannot be prolonged
to a larger open subscheme, such that the line bundle $\sigma^*(\OO_F(-1))$
on $U$ admits a prolongation to a line bundle on $X$. 
\end{proposition}

{\it Proof.} 
For any open $U\subset X$, 
we have a natural bijection between the set of all 
sections $\sigma: U \to {\bf P}(F)$ and the set of 
all line subbundles $f' : L' \subset F|U$, given by pulling back the 
tautological line subbundle $\OO_F(-1)\hra F_{{\bf P}(F)}$ over ${\bf P}(F)$.
This defines a natural map $R(F/X/S) \to \Sigma(F/X/S)$. 
If $L'$ can be prolonged to a line bundle $L$ on $X$, then  
by Lemma \ref{locally free sheaves are Z closed} below,  
any such a prolongation 
(that is, a pair $(L, L'\stackrel{\sim}{\to} L|U)$) is unique 
up to a unique isomorphism and moreover if such an $L$ exists then   
$f'$ admits a unique prolongation to a homomorphism 
$f: L\to F$. This shows that the natural map $R(F/X/S) \to \Sigma(F/X/S)$
is a bijection. \hfill$\square$

\begin{lemma}\label{locally free sheaves are Z closed} 
Let $\pi: X\to S$ be a smooth morphism where
$S$ is noetherian, let $j: U\hra X$ be an open subscheme
which is relatively big over $S$, and
let ${\mathcal E}$ be a locally free $\OO_X$-module. Then
the homomorphism ${\mathcal E} \to j_*({\mathcal E}|U)$ is an isomorphism.
\end{lemma}

{\it Proof.} For any $z\in Z = X-U$, if $\pi(z) = s \in S$, then 
$\depth(\OO_{X_s,z}) \ge 2$ as $\OO_{X_s,z}$ is a regular local ring of dimension 
$\ge 2$. By EGA IV$_2$ Proposition 6.3.1, we have
$\depth(\OO_{X,z}) = \depth(\OO_{S,s})+ \depth(\OO_{X_s,z})$, hence
$\depth(\OO_{X,z}) \ge 2$. Therefore, $\depth_Z(\OO_X) = \inf_{z\in Z} 
\depth(\OO_{X,z}) \ge 2$. Hence the desired conclusion
follows from EGA IV$_2$ Theorem 5.10.5. \hfill$\square$

\rem\label{effect of tensoring by a line bundle} 
If $F$ is replaced by $F\otimes K$
for a line bundle $K$ on $X$, then note that we have a natural isomorphism
${\bf P}(F) \cong {\bf P}(F\otimes K)$ and a compatible 
natural bijection $R(F/X) \cong R((F\otimes K)/X)$
(or $R(F/X/S) \cong R((F\otimes K)/X/S)$ in the relative case)
between the sets of $F$-linearized and $F\otimes K$-linearized rational
sections of the projective bundles. 
These natural isomorphisms and natural bijections satisfy the $1$-cocycle
condition as $K$ varies.

\medskip

For any vector bundle $F$ on a noetherian integral locally factorial 
scheme $X$, we can consider the following three sets.

(1) The set $\Sigma_{\eta}$ of all {\bf generic sections} 
$\spec \kappa(\eta) \to {\bf P}(F)$ 
of $\pi : {\bf P}(F)\to X$ over
the generic point $\eta \in X$.

(2) The set $\Sigma$ of all {\bf maximal rational sections}, that is, all
pairs $(U, \sigma)$ where $U\subset X$ is a big 
open subscheme of $X$ and $\sigma : U \to {\bf P}(F)$ is a section 
of the projection $\pi : {\bf P}(F)\to X$ over $U$, such that 
$(U, \sigma)$ is maximal in the sense that $\sigma$ does not admit
a prolongation to a section $\sigma'$ of $\pi$ which is 
defined over a strictly larger open subscheme of $X$ containing $U$.

(3) The set $R(F/X)$ of all $F$-linearized rational 
sections of $\pi : {\bf P}(F)\to X$, defined above.

\begin{proposition}\label{rational sections}
For any vector bundle $F$ on a noetherian integral locally factorial 
scheme $X$, we have a
natural bijection $\Sigma \to \Sigma_{\eta}$ which
sends $(U, \sigma)$ to the value of $\sigma$ at $\eta$.
Moreover, we have a natural bijection 
from $R(F/X)$ to $\Sigma$ which sends $(L,f)$ to the pair $(U,\sigma)$,
which is well-defined and maximal hence belongs the set $\Sigma$, 
where $U =  \{ x\in X \,|\, \rank(f_x) = 1\}$, 
and $\sigma : U \to {\bf P}(F)$ 
is defined by the line subbundle $f|U : L|U \to F|U$. 
\end{proposition}

{\it Proof.} We will just comment that the assumption of 
local factoriality of $X$ allows a unique (up to unique isomorphism) 
extension of a line bundle from a big open subset $U$ to all of $X$.
Moreover, for a line bundle on $X$, a section over a big $U$ will uniquely
extend to a global section because of the same assumption implies normality. 
The rest of the proof is a routine exercise and we omit the details. 
\hfill $\square$

\bigskip

\stm\label{due to Najmuddin} 
{\bf Example of a non-prolongable line bundle.} (Due to Najmuddin Fakhruddin.) 
Let $k$ be any field, let $S = \spec k [\epsilon]/(\epsilon^2)$, 
and let $X = {\mathbb P}^2_S$.
Let  $X_k = {\mathbb P}^2_k$ and let $P_0\in X_k$ be a closed $k$-rational 
point. Let $U = X - \{ P_0\}$, which is an open subscheme relatively big 
over $S$. The following is an example of a 
line bundle on $U$ which does not prolong to a line bundle $L$ on $X$. 
A \v{C}ech calculation using an open cover of $U_k= X_k - \{ P_0\}$ by
two copies of ${\mathbb A}^2_k$ shows that $H^1(U_k,\OO_{U_k})\ne 0$.
By elementary first order deformation theory (see e.g. [Ni 2]
for an expository account), $H^1(U_k,\OO_{U_k})$
is the group of isomorphism classes of all line bundles on 
$U= U_k\times_k \spec\, k[\epsilon]/(\epsilon^2)$
whose restriction to $U_k$ is trivial. 
But as $H^1(X_k, \OO_{X_k}) =0$, 
any line bundle on $X= X_k\times_k \spec\, k[\epsilon]/(\epsilon^2)$ 
which is trivial on $X_k$ must be trivial. 
Hence any non-zero element of $H^1(U_k,\OO_{U_k})$ defines a line bundle on
$U$ 
which does not prolong to $X$. More explicitly, 
$Pic(U) = \Z\cdot [\OO_{X/S}(1)|U] \oplus H^1(U_k,\OO_{U_k})
= Pic(X)\oplus H^1(U_k,\OO_{U_k})$,
where the second summand is non-zero.

\section{Relative rational parabolic reductions}

\subsection{Some basics about reductive groups}

From now onwards, we fix a base field $k$ of characteristic zero, 
and a connected split reductive group scheme $G$ over $k$, along with 
a chosen maximal torus and Borel $T\subset B\subset G$,
where $T$ is split over $k$. 
We denote the corresponding set of 
simple roots by $\Delta \subset X^*(T)$, and
for each $\alpha \in \Delta$ we denote by $\omega_{\alpha} \in
\Q\otimes X^*(T)$ the 
corresponding fundamental dominant weight. Recall that the 
standard parabolic subgroups
$P$ of $G$ (that is, $P\supset B$) are in a one-one correspondence with 
subsets $I_P\subset \Delta$, where $I_P$ is the set of inverted roots
for $P$.

\stm\label{internal direct sum} 
Recall that there exist natural internal direct sum decompositions
\begin{eqnarray*}
\Q\otimes X^*(T) & = & (\Q\otimes X^*(P)|_T) \oplus 
(\oplus_{\alpha \in I_P} \Q \alpha), \mbox{ where, in turn,} \\ 
\Q\otimes X^*(P)|_T & = & (\oplus_{\alpha \in \Delta - I_P}
\Q\omega_{\alpha}) \oplus (\Q\otimes X^*(G)|_T). 
\end{eqnarray*}
In the above, $X^*(P)|_T$ and $X^*(G)|_T$ denote 
the images of the homomorphisms $X^*(P) \to X^*(T)$ and 
$X^*(G) \to X^*(T)$ on the character groups which are 
induced by the inclusions $T\hra P$ and $T\hra G$. 

For simplicity of notation,
we will denote the base-change of $G$ (or $T$ or $B$ etc.)
to any $k$-scheme $S$ (in particular, to an extension field $K$)
just by $G$ alone rather than by $G_S$ or $G_K$, 
when there is no danger of confusion. Note that for any field extension $K/k$, 
we have a natural identification $X^*(T_K) = X^*(T)$ 
(respectively, $X_*(T_K) = X_*(T)$) on the groups of characters 
(respectively, on the groups of $1$-parameter subgroups). Moreover,
note that the set of all parabolics $P\supset B$ in $G$ 
is in a natural bijection 
with the set of all parabolics $P'\supset B_K$ under $P\mapsto P_K$,
and both these sets are in a natural bijection with the set of all subsets of 
$\Delta$. Hence for notational simplicity, we will identify  
any standard parabolic in $G_K$ with the corresponding standard 
parabolic in $G$. 

Let $P$ be a standard parabolic. 
Then $G/P$ is projective, and 
the anti-canonical line bundle 
$\omega_{G/P}^{-1}$ of $G/P$ is a very ample line bundle which has a natural
$G$-action which lifts the $G$-action on $G/P$. The canonical line bundle
$\omega_{G/P}$ is the associated line bundle for the principal $P$-bundle
$G\to G/P$ for a multiplicative character which we denote by
$\lambda_P : P \to \GG_m$.
The corresponding weight $\lambda_P \in X^*(T)$, which lies in the 
negative ample cone of $G/P$, has the form 
$$\lambda_P = c\cdot \sum_{\alpha \in \Delta - I_P} \omega_{\alpha}$$
for some $c>0$. Then $V_P = \Gamma(G/P, \omega_{G/P}^{-1})^{\vee}$ 
is an irreducible $G$-representation with highest weight
$\lambda_P$. For any non-zero weight vector $v\in V_P$ with weight 
$\lambda_P$ (such a $v$ is unique up to scalar multiple), 
the isotropy subgroup scheme
$I_{[v]}\subset G$ at the point $[v] \in {\bf P}(V_P)$ 
for the $G$-action on ${\bf P}(V)$
is $P$. This gives a $G$-equivariant closed imbedding 
$G/P \hra {\bf P}(V_P)$, under which $eP \mapsto [v]$.  
This is just the projective embedding given by the very ample 
line bundle $\omega_{G/P}^{-1}$. 

More generally, for any $\lambda: P \to \GG_m$ which as an element of $X^*(T)$ 
is a strictly positive linear combination of the $\omega_{\alpha}$ for
$\alpha \in \Delta - I_P$, 
(that is, $\lambda$ is in the negative ample cone of $G/P$), 
we similarly get a $G$-equivariant closed imbedding 
$G/P \hra {\bf P}(V_{\lambda})$ where $V_{\lambda}$ is an 
irreducible $G$-representation with highest weight
$\lambda$.

{}

Recall that the {\bf standard partial order} on $\Q\otimes X_*(T)$ is defined 
as follows. 
If $\mu, \nu \in \Q\otimes X_*(T)$, then we say that $\mu \le \nu$ if 
$\langle \chi, \mu \rangle = \langle \chi, \nu \rangle$ for all 
$\chi \in X^*(G)$ and 
$\langle \omega_{\alpha} , \mu \rangle \le 
\langle \omega_{\alpha} , \nu \rangle $
for all simple roots $\alpha \in \Delta$, where 
$\omega_{\alpha} \in \Q\otimes X^*(T)$ denotes the fundamental 
dominant weight corresponding to $\alpha$,
and $\langle ~,~\rangle$ denotes the duality pairing. 
The closed positive Weyl chamber $\ov{C}\subset \Q\otimes X_*(T)$ is 
the subset consisting of all $\mu\in \Q\otimes X_*(T)$ such that 
$\langle \alpha, \mu\rangle \ge 0$ for all $\alpha \in \Delta$.
The above partial order on $\Q\otimes X_*(T)$ induces the 
{\bf standard partial order} on $\ov{C}$.

\bigskip

\rem\label{positive expansion} 
If $I\subset \Delta$ is any subset, then we have the following elementary 
facts.\\
(i) Any $\alpha \in \Delta - I$ is expressible as  
$\alpha = \sum_{\beta \in I} b_{\beta} \beta  
+ \sum_{\gamma \in \Delta -I} c_{\gamma} \omega_{ \gamma}$ where each coefficient
$b_{\beta}$ is $\le 0$, and \\
(ii) any fundamental weight $\omega_{\alpha}$ 
is expressible as 
$\omega_{\alpha} = \sum_{\beta \in I} b_{\beta} \beta  
+   \sum_{\gamma \in \Delta -I} c_{\gamma} \omega_{ \gamma}$
where all coefficients $b_{\beta}$ and $c_{\gamma}$ are $\ge 0$. \\
The statement (i) is the statement labelled (R) in [Gu-Ni 1].  
The statement (ii) is immediate from Lemma 1.1 of Raghunathan [Rag],
and will be used below.

{}

\begin{lemma}\label{roots weights inequality}
Let $P$ be a standard parabolic in $G$. 
Let $\mu, \nu \in \ov{C}$, such that 
$\langle \chi, \mu\rangle = \langle \chi,\nu\rangle = 0$ for
all $\chi \in I_P$, and 
$\langle \chi, \mu\rangle = \langle \chi,\nu\rangle$ for
all $\chi \in \wh{G}$. Then $\mu \le \nu$ if 
and only if for all dominant weights 
$\lambda \in X^*(T)$ 
which lie in the negative ample cone of $G/P$, we have
$\langle \lambda ,\, \mu\rangle 
\le \langle \lambda ,\, \nu \rangle$.
\end{lemma}

{\it Proof.} 
By definition, for any $\mu, \nu \in \Q\otimes X_*(T)$, 
we have $\mu \le \nu$ if 
$\langle \chi, \mu \rangle = \langle \chi, \nu \rangle$ for all 
$\chi \in X^*(G)$ and 
$\langle \omega_{\alpha} , \mu \rangle \le 
\langle \omega_{\alpha} , \nu \rangle $
for all simple roots $\alpha \in \Delta$, where 
$\omega_{\alpha} \in \Q\otimes X^*(T)$ denotes the fundamental 
dominant weight corresponding to $\alpha$.  
By Remark \ref{positive expansion}, 
for any $\alpha \in I_P$, the corresponding
fundamental dominant weight $\omega_{\alpha}$ can be
written as a linear combination 
$\omega_{\alpha} = \sum_{\beta \in I_P} b_{\beta} \beta  
+   \sum_{\gamma \in \Delta -I_P} c_{\gamma} \omega_{ \gamma}$
where each $c_{\gamma} \in \Q^{\ge 0}$. 
As $\sum_{\beta \in I_P} b_{\beta} \beta$ evaluates to $0$ on both 
$\mu$ and $\nu$, we see that 
$\langle \omega_{\alpha},  \mu\rangle = 
\sum_{\gamma \in \Delta -I_P} c_{\gamma} \langle \omega_{ \gamma},  
\mu\rangle$
and $\langle \omega_{\alpha}, \nu\rangle 
= \sum_{\gamma \in \Delta -I_P} c_{\gamma} \langle \omega_{ \gamma},
\nu\rangle$, where each $c_{\gamma}\ge 0$.
Hence to prove that $\mu \le \nu$
it is necessary and sufficient to prove the following:

{\bf (*) } If $\alpha \in \Delta-I_P$ then  
$\langle  \omega_{\alpha},\, \mu \rangle 
\le \langle \omega_{\alpha} ,\, \nu \rangle$.

For each $\alpha \in \Delta-I_P$ and each integer $N \ge 0$,
we form an element 
$$\lambda_{\alpha, N} = N n_0\,\omega_{\alpha} + n_0\sum_{\beta \in 
\Delta - I_P} \omega_{\beta}$$
where $n_0 \in \Z^{>0}$ is so chosen that  
$n_0 \omega_{\alpha} \in X^*(T)$ for all
$\alpha \in \Delta$.
By taking limit as $N\to \infty$, it follows that 
in order to prove {\bf (*)}, it is enough to
prove the following:

{\bf (**) } If $\alpha \in \Delta-I_P$ and $N\ge 0$, then 
$\langle \lambda_{\alpha, N} ,\, \mu\rangle 
\le \langle \lambda_{\alpha, N} ,\, \nu \rangle$.
 
Note that each such element $\lambda_{\alpha, N}$ 
is a dominant weight in $X^*(T)$, 
which lies in the negative ample cone of $G/P$.
Hence by hypothesis, the inequality  {\bf (**) } holds.
\hfill$\square$

\subsection{Principal bundles and linearized rational $P$-reductions}

If $E\to X$ is a principal $G$-bundle on a noetherian $k$-scheme $X$, then for 
any standard parabolic $P$,  
the $G$-equivariant closed imbedding $G/P \hra {\bf P}(V_P)$ 
defined above gives rise to a closed embedding
$$i : E/P \hra {\bf P}(E(V_P))$$
where $E(V_P)$ is the vector bundle on $X$ associated to $E$ by the 
$G$-representation $V_P$, and ${\bf P}(E(V_P))$ is the projective bundle 
of lines in $E(V_P)$. 
Alternately, we can regard ${\bf P}(E(V_P))$ as the
associated bundle of $E$ for the $G$-action on the projective space 
${\bf P}(V_P)$. More generally, for any dominant weight 
$\lambda$ is in the negative ample cone of $G/P$, we get 
a closed embedding $E/P \hra {\bf P}(E(V_{\lambda}))$ where
$V_{\lambda}$ is an irreducible $G$-representation with highest weight $\lambda$.

To the above data, we now associate a subset
$$R(E/X, P) \subset R(E(V_P)/X)$$
which consists of all linearized rational sections
$[L, f: L \to E(V_P)]$ of ${\bf P}(E(V_P))$ over $X$ 
such that the section $\sigma$ of 
${\bf P}(E(V_P))$ defined by $f$ over the big open subscheme
$U = \{ x\in X\,|\, \rank(f_x) =1\}\subset X$ factors via the 
closed embedding $E/P \hra {\bf P}(E(V_P))$
defined above. We call any element $[L,f] \in R(E/X, P)$
as a {\bf linearized rational $P$-reduction} of 
the structure group of $E$.
More generally, for any dominant weight 
$\lambda$ is in the negative ample cone of $G/P$, we define
a subset $R(E/X, P, \lambda) \subset R(E(V_{\lambda})/X)$ 
by replacing $\lambda_P$ by $\lambda$ in the above. We
call its elements $(L, f: L\to E(V_{\lambda})$ as  
{\bf $\lambda$-linearized rational $P$-reductions} of $E$. 

The following proposition is an immediate consequence
of Proposition \ref{rational sections}.

\begin{proposition}\label{rational parabolic reductions}
Let $X$ be a noetherian locally factorial integral scheme
over $k$, let $E$ be a principal $G$-bundle on $X$, and 
$P$ be a standard parabolic in $G$.
Then the restrictions of the natural bijections of Proposition 
\ref{rational sections}
gives natural bijections \\
$\Sigma(E/X,P) \to \Sigma_{\eta}(E/X,P)$ 
and $R(E/X,P) \to \Sigma(E/X,P)$ 
between the following
three sets.

(1) The set $\Sigma_{\eta}(E/X,P)$ of all {\bf generic $P$-reductions}, 
that is, all sections 
$\spec \kappa(\eta) \to E/P$ 
of $\pi: E/P\to X$ over
the generic point $\eta \in X$.

(2) The set $\Sigma(E/X,P)$ of all 
{\bf rational $P$-reductions}, 
that is, all pairs $(P, \sigma)$ where 
$\sigma : U \to E/P$ is a section 
of the projection $\pi: E/P \to X$ over a big 
open subscheme $U\subset X$, such that 
$(P, \sigma)$ is maximal in the sense that $\sigma$ does not admit
a prolongation to a section $\sigma'$ of $\pi$ which is 
defined over a strictly larger open subscheme of $X$ containing $U$.

(3) The set $R(E/X,P)$ of all linearized
rational $P$-reductions defined above.

More generally, for any dominant weight 
$\lambda$ in in the negative ample cone of $G/P$, restriction 
gives a bijection $R(E/X,P, \lambda) \to \Sigma(E/X,P)$.
\hfill $\square$
\end{proposition}

{}

We now come to the relative case, where we make the following definition
which is inspired by the Proposition \ref{rational parabolic reductions}.

\defin\label{relative rational reduction} 
Let $X\to S$ be a morphism of noetherian $k$-schemes, and let 
$E$ be a principal $G$-bundle on $X$. 
A {\bf relative linearized rational $P$-reduction
of $E$ w.r.t. $S$} is a rational linearized $P$-reduction
$[L,f]$ of $E$ over $X$ such that 
the open subscheme $U= \{ x\in X\,|\,\rank(f_x)=1\} \subset X$
is relatively big over $S$. We denote by 
$$R(E/X/S, P)$$ 
the set of all such $[L,f]$. In terms of the Definition 
\ref{groupoid is equivalent to set} above, we have 
$$R(E/X/S, P) = R(E/X,P) \cap R(E(V_P)/X/S)
\subset R(E(V_P)/X).$$

\section{Canonical reductions and their restrictions}

Let $(X,\OO_X(1))$ be a connected smooth projective variety over 
a field $K$, together with a very ample line bundle.
Let $U\subset X$ be a big open subscheme, and let $F$ be a principal
$H$-bundle on $U$, where $H$ is a reductive group scheme over $K$. 
Recall that $F$ is said to be {\bf semi-stable} if for any
parabolic $Q\subset H$, any section $\sigma: W\to F/Q$ defined on
a big open subscheme $W$ of $U$, and any dominant character
$\chi : Q \to \GG_{m,K}$, we have 
$$\deg(\chi_*\sigma^*F) \le 0$$
where $\sigma^*E$ is the principal $P$-bundle on $W$ 
defined by the reduction $\sigma$,
and $\chi_*\sigma^*E$ is the $\GG_m$-bundle obtained by
extending its structure group via $\chi: Q \to \GG_m$,
which is equivalent to a line bundle on $W$. This line
bundle extends uniquely (up to a unique isomorphism) to a line bundle
on $X$, denoted again by $\chi_*\sigma^*E$, and 
$\deg(\chi_*\sigma^*F)$ is its degree w.r.t. $\OO_X(1)$.

Now let $K$ be an extension field of $k$, let
$(X,\OO_X(1))$ be a connected smooth projective variety over $K$,
and let $E$ be a principal $G$-bundle on $X$. 
Recall that a {\bf canonical reduction} of $E$ is  
a rational $P$-reduction $(P,\sigma)$
of structure group of $E$ to a standard parabolic $P\subset G$
(see Proposition \ref{rational parabolic reductions}.(2))
for which the following two conditions (C-1) and (C-2) hold:

{\bf (C-1)} If $\rho: P \to L = P/R_u(P)$ is the Levi quotient of $P$ 
(where $R_u(P)$ is the unipotent radical of $P$) then the principal $L$-bundle 
$\rho_*\sigma^*E$ is a semistable principal $L$-bundle defined 
on the big open subscheme $U$ on which $\sigma$ is defined. 

{\bf (C-2)} For any non-trivial character $\chi: P \to \GG_m$ 
whose restriction to the chosen maximal torus $T\subset B \subset P$ 
is a non-negative linear combination $\sum n_i\alpha_i$ 
of simple roots $\alpha_i \in \Delta$ (where $n_i \ge 0$, and at least
one $n_i \ne 0$), we have $\deg(\chi_*\sigma^*E) > 0$.

Recall that the {\bf type} of any  rational $P$-reduction 
$\sigma: U\to E/P$ of $E$ is  
the element $\mu_{(P,\sigma)}(E)\in \Q\otimes X_*(T)$ defined 
w.r.t. decomposition given in statement \ref{internal direct sum} by 
\begin{eqnarray*}
\langle \chi, \mu_{(P,\sigma)}(E)\rangle & = & \left\{
\begin{array}{ll}
\deg(\chi_*\sigma^*E) & \mbox{if }  \chi \in X^*(P)|_T ,  \\
0                     & \mbox{if } \chi \in I_P.
\end{array}\right.
\end{eqnarray*}

In particular, if $(P,\sigma)$ is a canonical reduction, then
the {\bf Harder-Narasimhan type} $\hn(E)$ is defined to be 
$\mu_{(P,\sigma)}(E)\in \Q\otimes X_*(T)$. This is well-defined, as
a canonical reduction is unique.
It can be shown that $\hn(E) \in \ov{C} \subset \Q\otimes X_*(T)$.

\stm\label{maximality property of the canonical reduction}
{\bf Maximality property of the canonical reduction.} 
Let $E$ be a principal $G$-bundle on a smooth connected projective 
variety $(X,\OO_X(1))$ and let 
$\sigma : U\to E/P$ define a rational 
reduction to a standard parabolic $P$. Then $\mu_{(P,\sigma)} \le \hn(E)$
in $\Q\otimes X_*(T)$, and if equality holds then
$\sigma : U\to E/P$ defines the canonical reduction.

\begin{lemma}\label{fundamental inequality}
{\bf (Fundamental Inequality.) }
Let $(X,\OO_X(1))$ be a smooth connected projective variety over an
extension field $K$ of $k$.
Let $E$ be a principal $G$-bundle on $X$, let $P$ be a standard
parabolic in $G$, and let $\sigma : U\to E/P$ be a section
over a big open subscheme $U\subset X$. Suppose that 
$U$ is the maximal open subscheme to which $\sigma$ can prolong.
Let $\mu_{(P,\sigma)}(E) \in \Q\otimes X_*(T)$ be the type of the reduction 
$\sigma : U\to E/P$. 
Let $Y\subset X$ be any smooth connected 
hypersurface such that $Y \in |\OO_X(m)|$ for some $m \ge 1$
(note that $m = \deg(Y)/deg(X)$). 
Let $U'\subset Y$ be a big open subscheme of $Y$ such that 
$U\cap Y \subset U'$, and let $\sigma' : U' \to (E|Y)/P$
be a section which prolongs $\sigma|(U\cap Y)$ 
(such a prolongation of $\sigma$ exists, as follows by
applying the valuative criterion for properness to 
$(E|Y)/P \to Y$). Let 
$\mu_{(P,\sigma')}(E|Y) \in \Q\otimes X_*(T)$ be the type of the reduction 
$\sigma' : U'\to (E|Y)/P$. Then we have
$$m\cdot \mu_{(P,\sigma)}(E) \le \mu_{(P,\sigma')}(E|Y).$$
Moreover, in the above we have equality 
$m\cdot \mu_{(P,\sigma)}(E) = \mu_{(P,\sigma')}(E|Y)$
if and only if $U\cap Y$ is big in $Y$.
\end{lemma}

{\it Proof.} 
By definition of the type of a rational parabolic reduction, we have 
$\langle \chi ,\, m\cdot\mu_{(P,\sigma)}(E)\rangle 
= \langle \chi ,\, \mu_{(P,\sigma')}(E|Y)\rangle =0$
for all $\chi \in I_P$. 

If $\chi\in \wh{G}$, then for any rational parabolic reduction
$(P,\sigma)$ of $E$ and $(P,\sigma')$ of $E|Y$, 
we have isomorphisms of line bundles $\chi_*(E) \cong \chi_*(\sigma^*(E))$
over $X$, and $\chi_*(E)|Y \cong \chi_*(E|Y) \cong \chi_*({\sigma'}^*(E|Y))$ 
over $Y$. Also, for any line bundle $L$ on $X$, we have 
$\deg_Y(L|Y) = m\cdot \deg_X(L)$. 
Hence 
$$\langle \chi , m\cdot \mu_{(P,\sigma)}(E)\rangle
= m \cdot \deg_X(\chi_*(E)) = \deg_Y(\chi_*(E|Y)) =   
\langle \chi ,\mu_{(P,\sigma')}(E|Y)\rangle.$$ 
Hence by Lemma \ref{roots weights inequality}, 
the inequality $m\cdot \mu_{(P,\sigma)}(E) \le \mu_{(P,\sigma')}(E|Y)$
will follow if we prove the following:

{\bf (***) } If $\lambda \in X^*(T)$ is a dominant weight 
which lies in the negative ample cone of $G/P$,
then $m \cdot \langle \lambda ,\, \mu_{(P,\sigma)}(E)\rangle 
\le \langle \lambda ,\, \mu_{(P,\sigma')}(E|Y)\rangle$.

To see this, let $L$ denote the line bundle 
$\lambda_*\sigma^*(E)$ on $X$. Under
the bijection given by Proposition \ref{rational parabolic
reductions}, the rational $P$-reduction $(P,\sigma)$ corresponds
to $[L,f] \in R(E/X, P, \lambda)$ where
$f : L = \lambda_*\sigma^*(E) \to E(V_{\lambda})$ is an injective
$\OO_X$-linear homomorphism, and  
$$U = \{ x\in X\,|\, \rank(f_x) = 1\}.$$

Let $Y\subset X$ be any smooth connected 
hypersurface such that $Y \in |\OO_X(m)|$ for some $m \ge 1$
(note that $m = \deg(Y)/\deg(X)$). 
Let $U'\subset Y$ be a big open subscheme of $Y$ such that 
$U\cap Y \subset U'$, and let $\sigma' : U' \to (E|Y)/P$
be a section with 
$$\sigma' |(U\cap Y) = \sigma| (U\cap Y).$$
Such a pair $(U', \sigma')$ indeed exists, as follows by
applying the valuative criterion for properness to 
$(E|Y)/P \to Y$. Moreover, we can take $U'\subset Y$ to be
the maximal open subscheme to which $\sigma' : U' \to (E|Y)/P$
prolongs. 

By Proposition \ref{rational parabolic reductions}, 
the rational $P$-reduction $(P,\sigma')$ of $E|Y$ corresponds
to $[L',f'] \in R((E|Y)/Y, P, \lambda)$ where
$f' : L' = \lambda_*{\sigma'}^*(E|Y) \to (E|Y)(V_{\lambda})$ is an injective
$\OO_Y$-linear homomorphism, and  
$$U' = \{ y\in Y\,|\, \rank(f'_y) = 1\}.$$

As $\sigma'$ prolongs $\sigma| (U\cap Y)$, it follows that the 
homomorphism $f|Y: L|Y \hra (E|Y)(V_{\lambda})$ factors uniquely 
via $f' : L' \hra (E|Y)(V_{\lambda})$ to give an injective
$\OO_Y$-linear homomorphism 
$$\varphi : L|Y \hra L'$$
which is fiberwise injective on $U\cap Y$. It now follows from
Remark \ref{degree inequality for line bundles} below that 
$$\deg(L|Y) \le \deg(L')$$
and equality holds if and only if the set
$Z = \{ y\in Y \,|\, \varphi_y =0 \}$ is of codimension $\ge 2$
in $Y$. As $f_y = f'_y\circ \varphi_y$, we have 
$Z = Y - (U\cap Y)$. Hence the equality $\deg(L|Y) = \deg(L')$
holds if and only if $U\cap Y$ is big in $Y$.

Now note that $\deg(L|Y) = m\cdot \deg(L)$, and
$\deg(L) = \langle \lambda ,\, \mu_{(P, \sigma)}(E)\rangle$ while 
$\deg(L') = \langle \lambda ,\mu_{(P, \sigma')}(E|Y)\rangle$.
Hence we get $m\cdot \langle \lambda ,\, \mu_{(P, \sigma)}(E)\rangle
\le \langle \lambda ,\mu_{(P, \sigma')}(E|Y)\rangle$.
This completes the proof of {\bf (***)}. 

Moreover, as $\deg(L|Y) = \deg(L')$ if and only if $U\cap Y$ is big
in $Y$, we see that  
$m\cdot \langle \lambda ,\, \mu_{(P, \sigma)}(E)\rangle
= \langle \lambda ,\mu_{(P, \sigma')}(E|Y)\rangle$
if and only if $U\cap Y$ is big in $Y$.
Hence, Lemma \ref{fundamental inequality} is proved. 
\hfill$\square$

{}

\rem\label{degree inequality for line bundles}
{\bf Degree inequality for line bundles.}
Let $(Y, \OO_Y(1))$ be an irreducible smooth projective variety over a 
field $K$, together with a very ample line bundle. If $\L_1\subset \L_2$
are coherent $\OO_Y$-modules such that both $\L_1$ and $\L_2$ are
invertible sheaves (line bundles), then 
$$\deg(\L_1)\le \deg(\L_2).$$ 
Moreover, equality holds in the above if and only if $\L_1 = \L_2$.
We leave the proof of this elementary fact to the reader.

{}

\begin{proposition}\label{HN-type rises under restriction}
{\bf HN-type rises under restriction.}
Let $(X,\OO_X(1))$ be a smooth connected projective variety over an
extension field $K$ of $k$, and let $Y\subset X$ be any smooth connected 
hypersurface such that $Y \in |\OO_X(m)|$ for some $m \ge 1$. 
Let $E$ be a principal $G$-bundle on $X$. Then we have
$$m\cdot \hn(E) \le \hn(E|Y).$$
Moreover, in the above we have equality 
$$m\cdot \hn(E) = \hn(E|Y)$$
only if for the canonical reduction 
$\sigma : U\to E/P$ of $E$, the intersection $U\cap Y$ is big in $Y$,
and $\sigma|(U\cap Y) : U\cap Y \to (E|Y)/P$ represents
the canonical reduction of $E|Y$.
\end{proposition}

{\it Proof.}
Let the section $\sigma : U\to E/P$ define a canonical reduction for $E$,
where $U\subset X$ is a big open subscheme. 
Moreover, assume that $U$ is the maximal open subscheme of $X$ 
to which $\sigma$ prolongs. 
Then using the notation of
Lemma \ref{fundamental inequality}, we have
\begin{eqnarray*}
m \cdot \hn(E) 
& = & m \cdot \mu_{(P,\sigma)}(E) 
            \mbox{ by definition of HN-type},\\
& \le & \mu_{(P, \sigma')}(E|Y) 
            \mbox{ by Lemma \ref{fundamental inequality}},\\
& \le & \hn(E|Y) \mbox{ by the maximality property
                 \ref{maximality property of the canonical reduction} 
                 of the HN-type}.
\end{eqnarray*}
By the last part of Lemma \ref{fundamental inequality}, 
equality holds in the above if and only if $U\cap Y$ is big in $Y$. 
Now by the maximality property of the canonical reduction
(see statement \ref{maximality property of the canonical reduction}),
it follows that $\sigma|(U\cap Y) : U\cap Y \to (E|Y)/P$ represents
the canonical reduction of $E|Y$.
This completes the proof of 
Proposition \ref{HN-type rises under restriction}. 
\hfill$\square$

\bigskip

By combining the analog for principal bundles of the 
semistable restriction theorem of Mehta and Ramanathan
with the Proposition \ref{HN-type rises under restriction} above,
we get the following result for HN-types of 
restrictions of principal bundles.

\begin{proposition}\label{Mehta-Ramanathan theorem for HN-types}
{\bf (Mehta-Ramanathan restriction theorem for HN-types.) }
Let $(X,\OO_X(1))$ be a smooth connected projective variety over
an extension field $K$ of $k$, of dimension $\ge 2$. 
Let $E$ be a principal $G$-bundle
on $X$. Then there exists an integer $m_0\ge 1$ such that for any 
general smooth hypersurface $Y\subset X$ where $Y \in |\OO_X(m)|$
where $m\ge m_0$, the HN-types $\hn(E)$ and $\hn(E|Y)$ 
are related by 
$$m\cdot \hn(E) = \hn(E|Y).$$
Moreover, if the canonical reduction of $E$ is 
the  rational reduction $(U,\sigma)$ to 
a standard parabolic $P$, then $U\cap Y$ is big in $Y$,
and the canonical reduction of $E|Y$ is represented by
the rational section of $(E|Y)/P$ defined by $\sigma|(U\cap Y)$.
\end{proposition}

\bigskip

\section{Restriction of relative canonical reduction}

Let $X\to S$ be a smooth projective morphism 
of noetherian $k$-schemes with geometrically connected fibers, 
with a given relatively very ample line bundle 
$\OO_{X/S}(1)$ on $X$. 

\defin\label{definition of relative canonical reduction}
 Let $E$ be a principal $G$-bundle on $X$, and suppose that 
the HN-type is constant over $S$, that is, there is some 
$\tau \in \ov{C}$ such that $\hn(E_s) = \tau$ for all $s\in S$.
A {\bf relative canonical reduction of $E$ w.r.t. $S$} will mean  
a relative linearized rational $P$-reduction 
$f: L\to E(V_P)$ of $E$ w.r.t. $S$ 
such that for each $s\in S$, the restriction 
$f_s: L_s\to E_s(V_P)$ is a canonical reduction of $E_s$ on 
$(X_s, \OO_{X_s}(1))$ over the base field $\kappa(s)$.

\begin{proposition}\label{restriction of relative canonical reduction} 
With $X\to S$ as above, let $E$ be a principal $G$-bundle on $X$
with constant HN-type $\hn(E_s) = \tau$ for all $s\in S$. 
Suppose that we are given a 
relative canonical reduction $(L,f)$ of $E$ w.r.t.  $S$, as defined above.  
Let $Y\subset X$ be a relative effective divisor such that 
$\OO_X(Y)\cong \OO_{X/S}(m)$, and such that $Y$ is smooth over $S$. 
Suppose that $\hn(E_s|Y_s) = m\cdot \tau$ for all $s\in S$. Then 
$(L|Y, f|Y)$, where 
$f|Y : L|Y \to E(V_P)|Y
= (E|Y)(V_P)$ denotes the restriction of $f$ to $Y$, 
is a relative canonical reduction of the principal $G$-bundle
$E|Y$ on $Y$ w.r.t. $S$. 
\end{proposition}

{\it Proof.}
Let $U\subset X$ be the subset where $\rank(f_x) =1$. 
Let $\sigma: U\to E/P$ be the section
defined by $f$. 
For any $s\in S$, we get a big open subset
$U_s = U\cap X_s$ and a section $\sigma_s = \sigma|U_s : U_s\to E_s/P$.
Note that $\sigma_s$ is defined by $f_s: L_s \to E_s(V_P)$. 
By Proposition \ref{rational parabolic reductions}, $(U_s, \sigma_s)$
is maximal. By assumption, $\hn(E_s|Y_s) = m\cdot \hn(E)$.
Hence by Proposition \ref{HN-type rises under restriction}, 
$U_s\cap Y_s$ is big in $Y_s$. Hence 
$f|Y : L|Y \to (E|Y)(V_P)$ 
is a relative canonical reduction of $E|Y$. \hfill$\square$

\rem In the light of Proposition 
\ref{restriction of relative canonical reduction}, 
if Theorem \ref{main theorem} is true,
then $S^{\tau}(E)$ must be a closed subscheme of $S^{m\tau}(E|Y)$.
The scheme $S^{m\tau}(E|Y)$ exists by inductive hypothesis 
on relative dimension. 
In what follows, we show that there indeed exists a closed
subscheme of $S^{m\tau}(E|Y)$ which has the universal property
required of $S^{\tau}(E)$.

\section{Embeddings of relative Picard schemes}

We will use repeatedly the following result of Mumford.

\begin{proposition}\label{use of m-regularity}
{\bf (Use of $m$-regularity.) } Let $S$ be any noetherian scheme, and let 
$\pi: X\to S$ be a projective morphism, together with a 
relatively very ample line bundle 
$\OO_{X/S}(1)$. Let $\F$ be any coherent sheaf on $X$. Then there exists
an integer $m$ such that for all $s\in S$, the sheaf 
$\F_s = \F|X_s$ is $m$-regular, that is, $H^i(X_s,\F_s(m-i)) =0$
for all $i\ge 1$.
\end{proposition}

{\it Proof (sketch).}
(See Mumford [Mu], or [Ni 1] for an expository account.)
There exists a surjection $\OO_X^n(-N) \to \F$, for some  
integers $n$ and $N$. Only finitely many different Hilbert
polynomials occur as Hilbert polynomials of $\F_s$ as $s$ varies over $S$.
Hence now the conclusion follows by Mumford's theorem 
(see page 101 of [Mu] or Theorem 5.3 of [Ni 1]). \hfill$\square$

\begin{proposition}\label{Picard embedding}
Let $X\to S$ be a smooth projective morphism of schemes 
of relative dimension $\ge 2$,
with geometrically connected fibers, 
where $S$ is noetherian, and let
$\OO_{X/S}(1)$ be a relatively very ample line bundle on $X/S$. 
Then there exists $m_0\in \Z$ with the following property:  
If $Y\subset X$ is a relative effective divisor
which is smooth over $S$ and such that 
$\OO_X(Y) \simeq \OO_{X/S}(m)$ where $m \ge m_0$, then the morphism
$$r: Pic_{X/S} \to Pic_{Y/S}$$
of relative Picard schemes which is induced by the inclusion $Y\hra X$  
is a closed embedding. 
\end{proposition}

{\it Proof.} 
We will assume the basics of Grothendieck's theory of relative Picard schemes, 
as explained in Kleiman [K]. 
For any integer $d$, let $Pic^d_{X/S}\subset Pic_{X/S}$ and 
$Pic^d_{Y/S}\subset Pic_{Y/S}$ denote the open and closed (`clopen')
subschemes corresponding to line bundles of degree $d$ 
w.r.t. the given relatively very ample line bundle $\OO_{X/S}(1)$ on $X/S$ 
and its restriction to $Y/S$. 
Note that $Pic_{X/S}$ and $Pic_{Y/S}$ are disjoint unions of the
subschemes $Pic^d_{X/S}$ and $Pic^d_{Y/S}$ respectively. 
The $S$-morphism $r: Pic_{X/S} \to Pic_{Y/S}$ maps $Pic^d_{X/S}$
into $Pic^{md}_{Y/S}$. 
Each restriction $r^d : Pic^d_{X/S} \to Pic^{md}_{Y/S}$ of $r$ 
is proper, as both $Pic^d_{X/S}$ and $Pic^{md}_{Y/S}$ are
projective over $S$ by the Kleiman boundedness theorem. 
Hence it follows that 
$r: Pic_{X/S} \to Pic_{Y/S}$ is proper. Moreover, if we show that each 
$r^d$ is a closed embedding, it would follow that so is $r$.

It is enough to prove the proposition after 
base change under an fppf (or surjective \'etale) morphism to $S$, 
as relative Picard schemes base change correctly under any morphism,
and a morphism is a closed embedding if and only if its
pullback under an fppf base change is a closed embedding.
Note that $Y\to S$ is an fppf morphism, and moreover, the base-change
$Y\times_SY \to Y$ of $Y\to S$ under $Y\to S$ 
admits a global section, namely, the diagonal. 
Hence without loss of generality, 
we can assume that $Y\to S$ admits a global section. 
Using this global section for normalization, we see that 
there exists a Poincar\'e line bundle $\LL$ 
on $X\times_S Pic^0_{X/S}$.

Let $\omega_{X/S} = \Omega^n_{X/S}$ 
be the rank $1$ locally free sheaf on $X$ of top relative exterior forms, 
where $n$ denotes the fiber dimension of $X/S$.
Let $p_1: X\times_S Pic^0_{X/S}\to X$ and 
$p_2: X\times_S Pic^0_{X/S}\to Pic^0_{X/S}$ denote the projections.
Consider the line bundles $\LL \otimes p_1^*\omega_{X/S}$
and $\LL^{-1} \otimes p_1^*\omega_{X/S}$ on $X\times_S Pic^0_{X/S}$.
As $Pic^0_{X/S}$ is noetherian and as 
$p_2: X\times_S Pic^0_{X/S} \to Pic^0_{X/S}$
is projective, by Mumford's theorem on $m$-regularity
(Proposition \ref{use of m-regularity})
applied to $\LL \otimes p_1^*\omega_{X/S}$
and $\LL^{-1} \otimes p_1^*\omega_{X/S}$, 
there exists an integer $m_1$ such that 
for all $m\ge m_1$,  
for any $z\in Pic^0_{X/S}$ and for all $i\ge 1$ we have
$$%
H^i(X_z, \LL_z^{-1}\otimes \omega_{X_z/z}(m)) = 
H^i(X_z, \LL_z\otimes \omega_{X_z/z}(m)) = 0,$$
where 
$X_z = X \times_S z \subset X\times_S Pic^0_{X/S}$
and $\LL_z = \LL|X_z$. 
Hence by Grothendieck duality on $X_z$ we have
$$%
H^{n-i}(X_z, \, \LL_z(-m)) = 
H^{n-i}(X_z, \, \LL_z^{-1}(-m)) = 0$$
for all $z\in Pic^0_{X/S}$, $m\ge m_1$ and $i\ge 1$.
Next by a similar argument for $X\to S$, we see that 
there exists an integer $m_2$ such that 
$$%
H^{n-i}(X_s, \, \OO_{X_s}(-m)) = 0
\mbox{ for all } s\in S, ~ m\ge m_2 \mbox{ and } i\ge 1.$$
We will show that the integer $m_0 = \max\{ m_1, m_2\}$ 
has the property required by the proposition.

By Lemma \ref{when embedding}, in order to prove the proposition it is
enough to prove that given any algebraically closed field $K$ with a 
morphism $\spec K \to S$, the base change 
$r^d_K: Pic^d_{X_K/K} \to Pic^{md}_{Y_K/K}$ is a closed embedding. 
As this morphism is a translate of $r^0_K$ by a $K$-point of 
$Pic^d_{X_K/K}$, it is enough to show that the proper morphism
$r^0_K: Pic^0_{X_K/K} \to Pic^0_{Y_K/K}$ is a closed embedding.

We first show that the abstract group homomorphism
$r^0_K(K): Pic^0_{X_K/K}(K) \to Pic^0_{Y_K/K}(K)$ is injective. 
As $K$ is algebraically closed, this is just the restriction map 
$Pic^0(X_K) \to Pic^0(Y_K)$. If a point $P \in Pic^0_{X_K/K}(K)$
is the isomorphism class of a line bundle $L$ on $X_K$,
then note that $L$ is isomorphic to the base change of 
the Poincar\'e bundle $\LL$ on $X\times_S Pic_{X/S}$
under the morphism $P: \spec K \to Pic_{X/S}$. 
Note that $Y_K$ is smooth projective over $K$,
and in fact $Y_K$ is connected as the relative dimension of 
$X/S$ is $\ge 2$ (see for example [H-AG] 3.7.9).
Hence if $L|Y_K$ is trivial, then 
$$\dim_K H^0(Y_K, L|Y_K) =  \dim_K H^0(Y_K, L^{-1}|Y_K) = 1.$$
But from our choice of $m_0$, for each $m\ge m_0$ we have
the cohomology vanishings 
{\footnotesize 
$H^0(X_K, L(-m)) = H^1(X_K, L(-m)) = 0$ and 
$H^0(X_K, L^{-1}(-m)) =  H^1(X_K, L^{-1}(-m)) = 0$.}
Hence from the long exact cohomology sequences for 
the short exact sequences 
$0\to L(-m) \to L\to L|Y_K \to 0$ and 
$0\to L^{-1}(-m) \to L^{-1}\to L^{-1}|Y_K \to 0$,
it follows that the restriction maps 
$H^0(X_K, L) \to H^0(Y_K, L|Y_K)$ and 
$H^0(X_K, L^{-1}) \to H^0(Y_K, L^{-1}|Y_K)$ are bijections.
This shows that 
$$\dim_K H^0(X_K, L) =  \dim_K H^0(X_K, L^{-1}) = 1,$$
hence $L$ is trivial, as $X_K$ is integral and proper over $K$. 

Hence $r^0_K: Pic^0_{X_K/K} \to Pic^0_{Y_K/K}$ is a proper 
homomorphism of finite type group schemes over $K = \ov{K}$,
which is injective on $K$-valued points. 
Hence its schematic 
kernel $N_K$ is a finite connected scheme, supported
over the identity $e = \spec K \subset Pic^0_{X_K/K}$. To show that 
$N_K$ is the identity, it is thus enough to show that the tangent map
$d_er^0_K : T_e Pic^0_{X_K/K} \to T_e Pic^0_{Y_K/K}$ is injective.
But this is just the restriction map 
$H^1(X_K, \OO_{X_K}) \to H^1(Y_K, \OO_{Y_K})$, 
which is injective for $m\ge m_0$ as its kernel 
$H^1(X_K, \OO_{X_K}(-m))$ is zero, so $N_K = \spec K$. 
This completes the proof that $r^0_K$ is a closed embedding,
and so the Proposition \ref{Picard embedding} is proved. \hfill$\square$

\begin{lemma}\label{when embedding} 
Let $S$ be a noetherian scheme, and let
$r: P\to Q$ be a proper morphism of $S$-schemes.
If for each $s\in S$,
the restriction $r_s: P_s\to Q_s$ is a closed embedding, then
$r$ is a closed embedding. Equivalently, if for each 
algebraically closed field $K$ with a morphism 
$\spec K \to S$, if the base-change $r_K: P_K\to Q_K$
is a closed embedding, then
$r$ is a closed embedding. 
\end{lemma}

{\it Proof.} 
As $r$ is proper, $r_*  \OO_P$ is coherent, and so the cokernel $C$ of 
$r^{\#} : \OO_Q \to r_*  \OO_P$ is a coherent sheaf of $\OO_Q$-modules. 
Being injective and proper, $r$ is finite hence affine. 
As $r_s$ is a closed embedding for each $s\in S$, we must have 
$C_s =0$ on $Q_s$ for each $s\in S$. 
Hence by the coherence of $C$, we get $C =0$, which shows that  
$r^{\#}$ is surjective. This completes the proof of the lemma.
\hfill$\square$

\rem\label{when lift of a line bundle is possible} 
Let $Y\subset X\to S$ be as before, so that we have a closed embedding
$r: Pic_{X/S}\hra Pic_{Y/S}$. For any 
line bundle $L_Y$ on $Y$, let $S_1\subset S$ be the closed
subscheme which is the schematic inverse image of 
$Pic_{X/S}\subset Pic_{Y/S}$ under the section $[L_Y] : S \to  Pic_{Y/S}$. 
It has the universal property that for any 
morphism $T\to S$, there exists a line bundle on 
$X_T$ which lifts (prolongs)
the line bundle $\phi^*(L_T)$ from $Y_T$ under the inclusion
$Y_T\hra X_T$ if and only if $T\to S$ factors via $S_1\hra S$.

\section{Proofs of the main results}

We will need the following two facts 
(Theorem \ref{cohomology vanishing} and Lemma \ref{direct image}) 
which are known to experts and can be found in more general 
forms in the literature 
(self-contained proofs are given in [Gu-Ni 2]).

\begin{theorem}\label{cohomology vanishing}
{\bf (Vanishing theorem for $R^i\pi_*F(-m)$ for $i\le n-1$.) }
Let $S$ be a noetherian scheme, let
$\pi: X\to S$ be smooth projective, and let 
$\OO_{X/S}(1)$ be a relatively ample line bundle on $X/S$.
Suppose that all fibers $X_s$ are connected, of 
constant dimension $n$. 
Then for any locally free sheaf $F$ on $X$, 
there exists an $m_0\in \Z$ such that 
for all $m\ge m_0$ and for all morphisms
$T\to S$ where $T$ is noetherian, we have 
$R^i(\pi_T)_*(F_T(-m)) =0$ for all $i\le n-1$,
where $F_T$ is the pullback of $F$ to $X_T = X\times_ST$
and $\pi_T : X_T \to T$ is the projection. \hfill$\square$
\end{theorem}

{}

\begin{lemma}\label{direct image} Let $S$ be a noetherian scheme and let 
$\pi: X\to S$ be a proper flat surjective morphism. Suppose moreover
that each fiber of $\pi$ is geometrically integral. Then the 
homomorphism $\pi^{\#}: \OO_S \to \pi_*\OO_X$ is an isomorphism.
\hfill$\square$
\end{lemma}

{}

\rem\label{what upper semicontinuity means} 
Let $S$ be a topological space, and let $({\mathcal P}, \le)$ be a 
partially ordered set. We will say that a map $h: S \to {\mathcal P}$ 
is {\bf upper semicontinuous} if 
for each $\tau \in {\mathcal P}$, the subset 
$\{ s\in S \,|\, h(s) \le \tau\}$ is open in $S$. 
It can be seen that a map $h: S\to {\mathcal P}$ is 
upper semicontinuous in the above sense if and only if for each 
$s_0 \in S$, there exists an open neighbourhood
$s_0 \subset U \subset S$ such that 
$h|U$ takes a maximum at $s_0\in U$, that is, 
$h(u) \le h(s_0)$ for all $u\in U$. Consequently, we have
$\ov{S^{\tau}} \subset \cup_{\tau' \ge \tau}\,S^{\tau'}$.

\bigskip

We now revert to the notation where $k$ is a given field of characteristic 
$0$, and $G$ is a given split reductive group scheme over $k$ 
together with a given split 
maximal torus and a Borel $T\subset B \subset X$.

\begin{proposition}\label{semicontinuity and uniqueness}
{\bf (Semicontinuity and uniqueness.) }
Let $X\to S$ be a smooth projective morphism of noetherian $k$-schemes
with geometrically connected fibers, together with a given
relatively very ample line bundle $\OO_{X/S}(1)$ on $X/S$. 
Let $E$ be a principal $G$-bundle on $X$.
Then we have the following.

(1) The function $S \to \ov{C} : s\mapsto \hn(E_s)$ is upper semicontinuous.
In particular, the semistable locus $S^{ss} = \{ s\in S\,|\,
\hn(E_s) =0\}$ is open in $S$.

(2) If the HN-type $\hn(E_s)$ is constant on $S$, 
then there exists at most one relative canonical reduction for $E$.
\end{proposition}

{\it Proof.}
The result is proved in [Gu-Ni 1] when $X/S$ has relative dimension $1$. 
We now proceed by induction on the relative dimension, which we 
assume is $\ge 2$. 
By application of $m$-regularity 
(see Proposition \ref{use of m-regularity}), there exists an integer
$m_1$ that each $\OO_{X_s}$ is $m_1$-regular
for all $s\in S$. 
Hence by Grothendieck's theorem on semicontinuity and base change
(see Theorem 5.10 in [Ni 1] for a formulation directly useful here), 
for all $m\ge m_1$, the direct image $\pi_*(\OO_{X/S}(m))$ is locally
free, the higher direct images $R^i \pi_*(\OO_{X/S}(m))$
are equal to $0$ for all $i\ge 1$, 
and for all $s\in S$, the evaluation homomorphism
$(\pi_*\OO_{X/S}(m))_s \to H^0(X_s, \OO_{X_s}(m))$
is surjective. 

For any $s\in S$ and a smooth effective divisor $H\subset X_s$
such that $H \in |\OO_{X_s}(m)|$, let $v\in H^0(X_s, \OO_{X_s}(m))$
define $H$. As $\pi_*(\OO_{X/S}(m))$ is locally
free and as $(\pi_*\OO_{X/S}(m))_s \to H^0(X_s, \OO_{X_s}(m))$
is surjective, there exists an open neighbourhood $U$ of $s$ in $S$
and a section $u\in (\pi_*\OO_{X/S}(m))(U)$ such that $u_s = v$. 
Let $Y\subset \pi^{-1}(U)$ be the subscheme defined by the vanishing of
$u$. As $v\in H^0(X_s, \OO_{X_s}(m))$ is a regular section, 
it follows by using the local criterion of flatness (see 
Lemma 9.3.4 of [K]) that $Y\to U$ is flat 
at all points of $Y_s$, hence by properness $Y\to U$ is flat 
over a neighbourhood of $s$. 
As $Y\to U$ is smooth over $s$, there exists a possibly 
smaller neighbourhood $s\in V \subset U$ such that $Y_V \to V$ is smooth.
Then $Y\to V$ is a relative effective divisor in $|\OO_{X_V/V}(m)|$, 
with $Y_s = H$. 
This shows how to prolong the smooth effective divisor $H\subset X_s$ to a 
smooth relative effective divisor in an open neighbourhood $V$ of $s$
in the linear system defined by $\OO_{X_V/V}(m)$.

{\it Proof of \ref{semicontinuity and uniqueness}(1)}:  
By Remark \ref{what upper semicontinuity means}, it is enough to 
prove that for any $s_0\in S$,  
the value $\tau = \hn(E_{s_0})$ is a local maximum for $\hn(E_s)$.
Given $s_0\in S$, by the Mehta-Ramanathan theorem 
(Proposition \ref{Mehta-Ramanathan theorem for HN-types}), 
there exists
an integer $m_0$ such that for any $m\ge m_0$, for a general 
smooth hypersurface $H \subset X_{s_0}$ in the linear system 
$|\OO_{X_0}(m)|$, the restriction $E_{s_0}|H$ has HN-type $m\cdot \tau$.
Let $m_{s_0} = \max\{m_0, m_1\}$ where $m_1$ was chosen
at the beginning of the proof. Let $m\ge  m_{s_0}$, and let
$H\subset X_{s_0}$ be a hypersurface in $|\OO_{X_0}(m)|$ 
which is connected and smooth over $\spec \kappa(s_0)$  
(which exists), 
with $\hn(E_{s_0}|H) = m\cdot \tau$. 

As we have shown above, after shrinking $S$ to a smaller
open neighbourhood of $s_0$ if needed, there exists a smooth 
relative effective divisor $Y\subset X$ over $S$ (where for simplicity,  
the same notation $S$ now denotes 
an open neighbourhood of $s_0$ in the original $S$), such that 
$\OO_X(Y) \cong \OO_{X/S}(m)$, and such that 
$$Y_{s_0} = H \subset X_{s_0}.$$ 
As $Y/S$ is of relative dimension less than that of $X/S$,
by inductive hypothesis on the relative dimension, 
the result holds for $E|Y$, hence the value
$m\cdot \tau$ at $s_0$ is a local maximum value for the 
function $s\mapsto \hn(E|Y_s)$. Hence by shrinking $S$ to a further
smaller open neighbourhood of $s_0$ 
(which we again denote just by $S$ for simplicity), we can assume that 
$m\cdot\tau$ is the maximum value of $\hn(E|Y_s)$ on $S$. 
By Proposition \ref{HN-type rises under restriction}, 
for any $s\in S$ we have $m\cdot \hn(E_s) \le \hn(E|Y_s)$.
Hence for all $s\in S$, we get the inequalities
$$m\cdot \hn(E_s) \le \hn(E|Y_s) \le \hn(E|H)
= m\cdot \hn(E_{s_0})$$ 
showing that $\hn(E_{s_0})$ 
is a local maximum for $\hn(E_s)$, as claimed.
This completes the proof of Proposition 
\ref{semicontinuity and uniqueness}.(1).

\medskip

{\it Proof of \ref{semicontinuity and uniqueness}.(2)}: 
Let $P$ be the unique standard parabolic in $G$ such that $\tau$ 
lies in the negative ample cone of $G/P$. 
We must show that if there exists two relative canonical
reductions $(L,f : L \to E(V_P))$ 
and $(L',f': L' \to E(V_P))$ for $E$, then there exists
an isomorphism $\phi: L\to L'$ such that $f = f'\circ \phi$. 

Note that such an isomorphism $\phi$, if it exists, is automatically unique
as $f'$ is an injective homomorphism of sheaves. Hence, it is enough to
prove the existence of $\phi$ 
after base change to an fppf cover of $S$,
as the resulting $\phi$ will satisfy the cocycle condition by uniqueness,
hence would descend to $S$. 
In particular, it is enough to prove the result
in a neighbourhood of 
each point $s_0\in S$.

By the proof of part (1) above, after possibly shrinking $S$ to a smaller
open neighbourhood of $s_0$ (which we again denote just by $S$), 
for any $m\ge m_{s_0}$ 
there exists a smooth relative effective divisor $Y\subset X$ over $S$, 
such that $\OO_X(Y) \cong \OO_{X/S}(m)$, and such that 
the inequalities    
$m\cdot \hn(E_s) \le \hn(E|Y_s) \le \hn(E|H) = m\cdot \hn(E_{s_0})$
hold for all $s\in S$, where $Y_{s_0} =H$ is chosen so as to satisfy 
the Mehta-Ramanathan theorem 
(Proposition \ref{Mehta-Ramanathan theorem for HN-types}). 
As by hypothesis $\hn(E_s)=\tau$ (constant) for all $s\in S$,
it follows that $\hn(E|Y_s) = m\cdot \tau$ for all $s\in S$. 
Moreover, by Proposition \ref{HN-type rises under restriction},
the restriction of the canonical reduction of $E_s$ to $Y_s$
gives the canonical reduction of $E_s|Y_s$ for each $s\in S$.
Hence if $(L,f : L \to E(V_P))$ 
and $(L',f': L' \to E(V_P))$ 
are two relative canonical
reductions for $E$, then $(L|Y,f|Y)$ 
and $(L'|Y,f'|Y)$ are two relative canonical
reductions for $E|Y$.
Hence by the inductive hypothesis, we must have an isomorphism 
$L|Y \cong L'|Y$ of line bundles on $Y$.

By Proposition \ref{Picard embedding}, 
the restriction morphism of the relative Picard schemes 
$$ r: Pic_{X/S} \to Pic_{Y/S}$$
is a closed embedding when $m$ is sufficiently large. 
Note that $L$ and $L'$ define global sections
$[L]$ and $[L']$ of $Pic_{X/S} \to S$, 
which map under the restriction morphism 
$r: Pic_{X/S} \to Pic_{Y/S}$ 
to the same global section $[L|Y] = [L'|Y]$ of
$Pic_{Y/S} \to S$. Hence $[L] = [L']$. 

Recall that we have made an fppf base change such that 
$Y\to S$ (and hence also $X\to S$) admits a global section. 
Hence from $[L] = [L']$ we conclude that 
there exists a line bundle $K$ on $S$ and an isomorphism 
$L\otimes \pi^*K \to L'$. By shrinking $S$ to a further smaller open
neighbourhood of $s_0$, we can assume that $K$ is trivial on $S$,
and hence we have an isomorphism $\psi : L\stackrel{\sim}{\to} L'$ on $X$. 
Hence we can replace $(L',f')$ by $(L, f'\circ \psi)$
as the second of the two relative canonical reductions of $E$ which we 
want to show to be isomorphic,
that is, we are now reduced to the case where $L= L'$, and
we have been given 
two relative canonical reductions $(L,f)$ and $(L, f')$ of $E$
with $[L|Y,f|Y]= [L|Y,f'|Y]$. Note that by the 
uniqueness of canonical reduction over $Y$, there
exists an automorphism $g: L|Y \to L|Y$ such that 
$f|Y = (f'|Y)\circ g$. Note that any automorphism $g$ 
of a line bundle on $Y$ is given by scalar 
multiplication by an invertible function $\gamma \in \Gamma(Y,\OO_Y)^{\times}$.
Hence we have $f|Y = \gamma \cdot(f'|Y)$.

As $\pi_Y = \pi|Y : Y\to S$ is smooth projective with geometrically 
connected fibers, in particular $\pi_Y$ is proper flat with geometrically
integral fibers. Hence by Lemma \ref{direct image}, the homomorphism
$\pi_Y^{\#}: \OO_S \to (\pi_Y)_*\OO_Y$ is an isomorphism of sheaves of rings.
It follows that any invertible function $\gamma \in \Gamma(Y,\OO_Y)^{\times}$
is of the form $h\circ \pi_Y$ where $h = \gamma\circ b  \in 
\Gamma(S,\OO_S)^{\times}$ where $b: S\to Y$ is a global section. 

By Theorem \ref{cohomology vanishing}, if $m$ is sufficiently
large, then as the relative dimension $n$ of $X/S$ is $\ge 1$, we 
have 
$$\pi_*(\hom(L, E(V_P)(-m))) = 0.$$
By tensoring with $\hom(L, E(V_P))$, 
the short exact sequence
$0 \to \OO_X(-m) \to \OO_X \to \OO_Y \to 0$ gives the short exact sequence
$$0 \to \hom(L, E(V_P))(-m) \to \hom(L, E(V_P))
\to \hom(L|Y, (E|Y)(V_P)) \to 0.$$
The associated long exact cohomology sequence shows that 
the restriction homomorphism 
$$\rho : H^0(X, \hom(L, E(V_P))) \to 
H^0(Y, \hom(L|Y, (E|Y)(V_P)))$$
is injective. As
$$\rho(f) = f|Y = h\cdot (f'|Y) = h\cdot \rho(f') = \rho(h\cdot f')$$ 
we conclude that $f = h\cdot f'$. Hence 
scalar multiplication by the element $h \in \Gamma(S,\OO_S)^{\times}$ 
gives an isomorphism 
$\phi : L\to L$ such that $f = f'\circ \phi : L\to E(V_P)$,
as desired. This completes the proof of Proposition 
\ref{semicontinuity and uniqueness}.(2).
\hfill$\square$

\rem\label{inverse image of cone}  
Let $F$ be a vector bundle on a scheme $X$, let
$Z \subset {\bf P}(F)$ be a closed subscheme, and let 
$f\in \Gamma(X, F)$ be a global section of $F$. 
Let $\pi: {\bf F} = \spec_{\OO_X}\sym_{\OO_X}(F^{\vee})\to X$ 
be the geometric vector bundle over $X$ associated to $F$, 
and let $f: X\to {\bf F}$
denote its section defined by $f$.  Let $\wh{Z} \subset {\bf F}$ 
denote the affine cone over $Z$. Let $W = f^{-1}(\wh{Z}) \subset X$
be its schematic inverse image under the section $f$. 
This has the property that if $\phi: X'\to X$ is any morphism of 
schemes, if $F' = \phi^*(F)$, if $Z' \subset {\bf P}(F')$
is the pullback of $Z$ and if $f'\in \Gamma(X', F')$ is the 
pullback of $f$, then 
$f'$ factors via $\wh{Z'}$
if and only if $f$ factors via $W\hra X$.

\begin{proposition}\label{lifting the relative canonical reduction}
{\bf (Lifting the relative canonical reduction)}
Let $X\to S$ be a smooth projective morphism of noetherian $k$-schemes
with geometrically connected fibers of dimension $\ge 2$, 
together with a given
relatively very ample line bundle $\OO_{X/S}(1)$ on $X/S$. 
Let $E$ be a principal $G$-bundle on $X$. 
Let $\tau$ be the global maximum of $\hn(E_s)$ over $S$.
Let $Y\subset X$ be a smooth relative effective divisor over $S$, 
such that $\OO_X(Y) \cong \OO_{X/S}(m)$, where $m$ is sufficiently large.
Let there exist a relative canonical reduction for the restriction  
$E|Y$ over $S$, with constant HN-type $m\tau$. Then there exists  
a unique closed subscheme $S'\subset S$ with the universal property that
for any morphism $T\to S$ where $T$ is noetherian, the pull-back
$E_T$ admits a relative canonical reduction of type $\tau$
if and only if $T\to S$ factors
via $S'$. Equivalently, we have \\
(i) there exists a relative canonical reduction 
for $E_{S'}$ over $S'$ of type $\tau$, and \\
(ii) if $T\to S$ is any noetherian base change such that $E_T$ admits
a relative canonical reduction of type $\tau$, then $T\to S$ factors
via $S'$.
\end{proposition}

{\it Proof.}
By uniqueness of relative canonical reductions
(Proposition \ref{semicontinuity and uniqueness}.(2)), 
a relative canonical reduction defined on an
fppf (or \'etale) cover of $S$ will descend to $S$, 
so it is enough to prove the result locally over an
fppf (or \'etale) cover of $S$. In particular,
we can assume that $S$ is connected.
Moreover, as explained in the course of the proof of 
Proposition \ref{semicontinuity and uniqueness}, 
we can assume that $Y$ has a global section over $S$.

By Proposition \ref{Picard embedding}, 
the restriction morphism of relative Picard 
schemes $r: Pic_{X/S} \to Pic_{Y/S}$ is a closed  
embedding if $m\gg 0$. So if $(L_Y,f_Y)$ is
a relative canonical reduction for $E|Y$, there exists a  
closed subscheme $S_1$ of $S$ 
on which $[L_Y]$ lifts to a section of $Pic_{X/S}$ and 
such that $S_1$ has the universal property  
described in Remark 
\ref{when lift of a line bundle is possible}. Let $X_1 = X\times_SS_1$,
$Y_1 = Y\times_SS_1$, and $E_1$ be the pullback of $E$ to $X_1$.
As $Y$ has a global section, so does $Y_1$.  
As $Y_1$ (and hence $X_1$) 
has a global section, this means $L_Y|Y_1$ lifts to a line bundle
$L_1$ over $X_1$, so we can assume $L_Y|Y_1 = L_1|Y_1$. 

For any $E_s$ on $X_s$ with $\hn(E_s)=\tau$,
let $d$ denote the degree of the line bundle $L_s$ for 
any canonical reduction $(L_s,f_s)$ of $E_s$. 
Then $d$ is constant as $S$ is connected (in fact, 
$d = \langle \lambda_P\,,\,\tau \rangle$).
Let $\LL$ denote the relative Picard line bundle 
on $X\times Pic_{X/S}^d$, normalized to be trivial on a chosen
global section of $X\to S$. 
Consider the projection $p_2: X\times Pic_{X/S}^d \to Pic_{X/S}^d$ and the 
vector bundle $F = \un{Hom}(\LL, p_1^*E(V_P))$ on 
$X\times Pic_{X/S}^d$. Applying 
Theorem \ref{cohomology vanishing} to the
above data gives an integer $m_0$ with the property 
given in the conclusion of Theorem \ref{cohomology vanishing}
for base changes $T\to Pic_{X/S}^d$. We base change to $T = S_1$,
under the morphism $[L_1] : S_1 \to Pic_{X/S}^d$. Note that
the pull-back
of $p_2 : X\times Pic_{X/S}^d \to Pic_{X/S}^d$ is 
$\pi_1: X_1\to S_1$, and Zariski locally over $S_1$, 
the pull-back of $F$ is 
isomorphic to $\un{Hom}(L_1, E_1(V_P))$. 
Hence by Theorem \ref{cohomology vanishing}
$R^{n-i}(\pi_1)_*\un{Hom}(L_1, E_1(V_P)(-m)) = 
0$ for $m \ge m_0$ and $i\ge 1$. As $n\ge 2$ by inductive hypothesis, 
both the $0$th and $1$st direct images vanish, and hence
the restriction map
$$\rho : H^0(X_1, \hom(L_1, E_1(V_P))) \to 
H^0(Y_1, \hom(L_1|Y_1, (E_1|Y_1)(V_P)))$$
is a bijection for $m \ge m_0$. Hence when $m$ is sufficiently large, 
there exists 
a unique $f_1 : L_1 \to E_1(V_P)$ such that $f_Y|Y_1 = f_1|Y_1$.
This gives a pair $(L_1,f_1: L_1\to E_1(V_P))$ 
which lifts (that is, prolongs) the pair 
$(L_Y,f_Y)$ from $Y_1$ to $X_1$. However, 
$(L_1,f_1)$ need not be a relative canonical reduction over $S_1$.
We will produce a maximal closed subscheme $S'\subset S_1$ over 
which it will be so.

By tensoring with $L_1^{-1}$, the homomorphism $f_1 : L_1 \to E_1(V_P)$ 
can be regarded as a section 
$f_1 \in \Gamma(X_1, E_1(V_P)\otimes L_1^{-1})$.
Note that we have a natural identification 
${\bf P}(E_1(V_P)) = {\bf P}(E_1(V_P)\otimes L_1^{-1})$
under which $E_1/P$ becomes a closed subscheme
$E_1/P \subset {\bf P}(E_1(V_P)\otimes L_1^{-1})$.
Let 
$W = f_1^{-1}(\wh{E_1/P}) \subset X_1$
be the closed subscheme of $X_1$ defined as in Remark 
\ref{inverse image of cone}.

Consider the composite $W \hra X_1 \to S_1$, which is a projective morphism. 
The flattening stratification of $S_1$ 
for this morphism 
has strata $S_{1, p(t)}$ indexed by the various Hilbert polynomials
$p(t)\in \Q[t]$ of the fibers w.r.t. $\OO_{X_1/S_1}(1)$. The
top stratum $S' = S_{1,p_0(t)}$, with Hilbert polynomial $p_0(t)$
equal to the 
Hilbert polynomial of the fibers of $X\to S$, is a closed
subscheme of $S_1$. We will show that $S'$ satisfies the proposition.

{\it Proof that $S'$ satisfies 
\ref{lifting the relative canonical reduction}.(i):}
To begin with, we will look at the points of $S'$. 
Let $(L',f')$ denote the base-change
of $(L_1,f_1)$ under $S'\hra S_1$.  
By definition, $s\in S'$ if and only if 
the fiber of $W\to S$ over $s$ is $X_s$, which means 
the section $f'_s \in \Gamma(X_s,E_s(V_P)\otimes L_s^{-1})$ 
lies in the affine cone $W_s$ over 
$E_s/P \subset {\bf P}(E_s(V_P)\otimes L_s^{-1})$.
Let $U\subset X_s$ be the open subscheme where $\rank(f'_x) =1$.
Then $U$ is non-empty and big in $X_s$,
as by assumption for each $s\in S'$, the intersection 
$U\cap Y_s$ is a big open subscheme of $Y_s$.
Hence $[L'_s, f'_s] \in R(E_s/X_s, P)$, that is,
$(L'_s, f'_s)$ is a linearized rational reduction of $E_s$. 

Let $\mu_{(P,\, \sigma)}(E_s)$ and $\mu_{(P,\, \sigma|(U\cap Y_s))}(E_s|Y_s)$ be the types
of the rational $P$-reductions  $\sigma : U \to E_s/P$ and 
$\sigma|(U\cap Y_s) : U\cap Y_s \to (E_s|Y_s)/P$ respectively. 
As $U$ and $U\cap Y_s$ are both big, by Lemma \ref{fundamental inequality}
we get the equality
$$m\cdot \mu_{(P,\, \sigma)}(E_s) = \mu_{(P,\, \sigma|(U\cap Y_s))}(E_s|Y_s).$$
As by assumption 
$\mu_{(P,\, \sigma|(U\cap Y_s))}(E_s|Y_s) = \hn(E_s|Y_s) = m\cdot \tau$,
it follows that 
$\mu_{(P,\, \sigma)}(E_s)\\  
= \tau$. Now by hypothesis of Proposition 
\ref{lifting the relative canonical reduction}, $\tau$ is the global 
maximum of $\hn(E_s)$ over $S$. Hence by the maximality property of 
canonical reductions (see Statement 
\ref{maximality property of the canonical reduction}), 
$\sigma : U \to E_s/P$ is the canonical reduction of 
$E_s$. Hence we have shown that $(L'_s, f'_s)$ defines the 
canonical reduction of $E_s$ at all points $s\in S'$. 

Let $X' = X\times_SS'$, $Y' = Y\times_SS'$ and $W' = W\times_SS'$. 
Note that by definition of $S'$ as the top stratum in the 
flattening stratification of $W\to S$ with Hilbert polynomial
that of fibers of $X\to S$, we see that $W'\subset X'$ is a closed
subscheme, and both $W' \to S'$ and $X'\to S'$ are flat projective 
with the same Hilbert polynomial $p_0$ of fibers. 
Hence we must have $W' = X'$. 
This shows that as schemes, $X'\subset W\subset X$. 
Let $E' = E|X'$. As above, $f': L' \to E'(V_P)$ 
denote the restriction of $f_1: L_1 \to E_1(V_P)$ to $X'$. 
As shown above, for each $s\in S'$ the restriction $f'_s = f'|X_s$  
has rank $1$ over a big open subscheme of $X_s$, which shows that 
$f': L' \to E'(V_P)$ has rank $1$ over a relatively big open
subscheme of $X'$ over $S'$. 
As $X'\subset W$, by Remark \ref{inverse image of cone}, the section 
of ${\bf P}(E'(V_P))$ 
defined by $(L',f')$ over the above relatively big open subscheme
factors via $E'/P \hra {\bf P}(E'(V_P))$. Hence 
$(L',f')$ is a linearized relative rational $P$-reduction 
of $E'$ over $S'$ in the sense of Definition 
\ref{relative rational reduction}. 
Now recall that we have shown that for each $s\in S'$, 
$(L',f')$ restricts to $X_s$ to give the canonical reduction of 
$E_s$. Hence by Definition \ref{definition of relative canonical reduction},
$(L',f')$ is a relative canonical of $E'$ over $S'$.
This completes the proof of the 
property (i) in the conclusion of Proposition 
\ref{lifting the relative canonical reduction}.

{\it Proof that $S'$ satisfies 
\ref{lifting the relative canonical reduction}.(ii):} 
Let $\phi: T\to S$ be a morphism, and let $(\L, g)$ be a relative canonical
reduction of $E_T$ on $X_T/T$.
Note that the pull-back of the given
relative canonical reduction $(L_Y,f_Y)$ of $E|Y$
under $\phi: T\to S$ is a relative
canonical reduction of $E_T|Y_T$ over $T$. 
Therefore by Propositions \ref{restriction of relative canonical reduction}
and \ref{semicontinuity and uniqueness} applied over $T$, 
the restriction of $(\L,g)$ to $Y_T$ is isomorphic to 
$(\phi^*L_Y,\phi^*f_Y)$.

In particular, as $\L$ is a prolongation of $\phi^*(L_Y)$ under 
$Y_T\subset X_T$, by 
Remark \ref{when lift of a line bundle is possible}, the morphism 
$\phi: T\to S$ factors via $S_1\hra S$, giving $\phi_1: T\to S_1$. 
The pullback of $f_1: L_1 \to E_1(V_P)$ under 
$\phi_1: T\to S_1$ gives another 
prolongation 
$(\phi_1^*L_1,\, \phi_1^*f_1)$
of $(\phi^*L_Y,\phi^*f_Y)$ under $Y_T\subset X_T$.
Zariski locally over $T$, we can identify $\L$ with $\phi_1^*L_1$
because of Proposition \ref{Picard embedding}. 
By Theorem \ref{cohomology vanishing},  
the restriction map
$$\rho : H^0(X_T, \hom(\L, E_T(V_P))) \to 
H^0(Y_T, \hom(\L|Y_T, (E_T|Y_T)(V_P)))$$
is injective,  
hence we must moreover have 
$(\phi_1^*L_1,\, \phi_1^*f_1)\cong (\L,g)$,
that is, there exists an isomorphism 
$u: \phi_1^*L_1 \to \L$ such that $\phi_1^*f_1 = g\circ u 
: \phi_1^*L_1 \to E_T(V_P)$.

As $(\L,g)$ is a relative canonical reduction, it follows 
that so is $(\phi_1^*L_1, \phi_1^*f_1)$. 
Hence $\phi_1^*f_1$ factors via the cone $\wh{E_T/P}$. 
Therefore 
by Remark \ref{inverse image of cone}, the morphism
$(\phi_1)_{X_1} : X_T \to X_1$ must factor through $W \hra X_1$. 
As the inclusion $W\hra X_1$ is monic, it follows
that $X_T\to T$ is the base-change of $W\to S_1$ under $\phi_1$.
As $X_T\to T$ is flat with Hilbert polynomial $p_0(t)$,
by the universal property of the flattening stratification for $W\to S_1$,
the morphism $\phi: T\to S_1$ factors via the 
corresponding stratum $S' = S_{1,p_0(t)}$
as we wished to show. This completes the proof of 
Proposition \ref{lifting the relative canonical reduction}.
\hfill $\square$

\bigskip

{\bf Proof of Theorem \ref{main theorem}.} 
With all the necessary ingredients in place at last, we can now complete
the proof of Theorem \ref{main theorem}.
We proceed by induction on the relative dimension of $X/S$.
We have earlier proved the theorem when $X/S$ is of relative dimension $1$
(see [Gu-Ni 1]). So now assume that the relative dimension is $\ge 2$.

It follows from the uniqueness of a relative canonical reduction
proved in Proposition \ref{semicontinuity and uniqueness}.(2)
that the assertion of the Theorem \ref{main theorem} is local over 
the base $S$, that is, it is enough to prove it 
in a neighbourhood of each point $s_0$ of $S$. As 
by Proposition \ref{semicontinuity and uniqueness}.(1) 
the HN-type $\hn(E_s)$ is upper semicontinuous on $S$, 
by shrinking the neighbourhood of $s_0$ if needed, we can assume that 
$\hn(E_s)$ attains a unique maximum value $\tau$ at $s_0$. 
Hence $|S|^{\tau}(E)$ is a closed subset of $|S|$. 
By the Mehta-Ramanathan theorem, given any sufficiently large $m\in \Z$, 
there will exist an effective divisor $H\subset X_{s_0}$
which is smooth over $s_0$ with $H \in |\OO_{X_{s_0}}(m)|$, 
such that the canonical reduction of $E_{s_0}$ on restriction to $H$
gives the canonical reduction of $E_{s_0}|H$.
This includes the property that the domain of definition
of the canonical reduction of $E_{s_0}$, which is a big open subset of $X_{s_0}$,
intersects $H$ to give a big open subset of $H$. 
 
If $m$ is taken to be sufficiently large, then after  
shrinking $S$ to a smaller neighbourhood of $s_0$ if needed, 
there will exist an effective relative divisor $Y\subset X$ over $S$
which is smooth over $S$ with $Y \in |\OO_{X/S}(m)|$, 
such that $Y_{s_0} = H$. 
In particular, we have $\hn(E|Y_{s_0}) = m\tau$.
Now by shrinking $S$ further, we can assume that for the 
restricted family $E|Y$ on $Y$ over $S$, 
$m\tau$ is the unique maximum value of $\hn(E|Y_s)$ on $S$. 

By our inductive hypothesis on the relative dimension of $X\to S$,
the Theorem \ref{main theorem} holds for the principal bundle $E|Y$ on $Y$ 
over $S$.
Hence there exists a closed subscheme $S^{m\tau}(E|Y) \subset S$ 
with the desired universal property for $E|Y$ for the type $m\tau$. 
For any $s\in |S|^{\tau}(E)\subset |S|$
we have $\hn(E|Y_s) \le m\tau$ as
$m\tau$ is the maximum for $\hn(E|Y_s)$ by assumption. 
On the other hand, by Proposition 
\ref{HN-type rises under restriction},  
if $s\in |S|^{\tau}(E)$ then the restriction of 
canonical reduction of $E_s$ to the subscheme $Y_s\subset X_s$ 
prolongs to a big open subset of $Y_s$ to give us a rational parabolic 
reduction of $E|Y_s$ of a type which is $\ge m\tau$. 
It follows that $\hn(E|Y_s) = m\tau$ whenever 
$s\in |S|^{\tau}(E)\subset |S|$ provided that $\tau$ is the unique 
maximum value of $\hn(E_s)$ and $m\tau$ is the unique 
maximum value of $\hn(E|Y_s)$ on $S$. Hence at the level of
sets, we have the inclusion $|S|^{\tau}(E)\subset |S|^{m\tau}(E|Y)$,
when $S$ is replaced by a neighbourhood of $s_0$ in $S$.

By Proposition \ref{restriction of relative canonical reduction}, 
with $S$ and $\tau$ as above, 
if $T\to S$ is a base change such that
$E_T$ admits a relative canonical reduction of constant type 
$\tau$, then this reduction restricts to 
$Y_T$ to give a relative canonical reduction of $E_T|Y_T$ of type $m\tau$. 
Hence $T\to S$ must factor through $S^{m\tau}(E|Y)\hra S$, 
so if the theorem is true, then $S^{\tau}(E)$ must be 
a closed subscheme of $S^{m\tau}(E|Y)$.

It therefore only remains to identify $S^{\tau}(E)$ as an 
appropriate closed subscheme of $S^{m\tau}(E|Y)$, and to show that it has
the desired universal property.
So we can replace the original $S$ by $S^{m\tau}(E|Y)$,
and assume (by the inductive hypothesis) that $E|Y$ 
has a relative canonical reduction
over $S$.
Moreover, we can assume (by replacing $S$ by an open neighbourhood
of $s_0$ if necessary)
that $\tau$ is the global maximum for $\hn(E_s)$ over $s$. 
The problem is to lift (prolong) this rational reduction of the 
structure group from $E|Y$ to $E$. 
By Proposition \ref{lifting the relative canonical reduction},
there exists a unique largest closed subscheme 
$S' \subset S^{m\tau}(E|Y)$ over which such a lift exists, and this  
closed scheme $S'$ has the functorial property which we would like
$S^{\tau}(E)$ to possess. This completes the proof of 
Theorem \ref{main theorem}. \hfill$\square$

\bigskip

The Theorem \ref{stacky main theorem}  can now be easily deduced from 
Theorem \ref{main theorem} 
as follows. 
Recall that the stack $Bun_{X/S}(G)$ of principal $G$-bundles on 
$X/S$ is the stack over $S$ whose objects over an $S$-scheme $T$ are all 
principal $G$-bundle $E$ on $X_T$. 
That $Bun_{X/S}(G)$ is an Artin stack can be seen as follows. 
The group $G$ can be embedded as a closed subgroup 
in $GL_{n,k}$ for some $n$, which gives a $1$-morphism 
$Bun_{X/S}(G)\to  Bun_{X/S}(GL_{n,k})$. By using an appropriate Hilbert scheme
of sections, it can be seen that the $1$-morphism  
$Bun_{X/S}(G)\to  Bun_{X/S}(GL_{n,k})$ is representable, in fact, schematic. 

For each $\tau \in \ov{C}$, we define an $S$-stack $Bun_{X/S}^{\tau}(G)$ 
whose objects over an $S$-scheme $T$ are all triples $(E, L, f)$ consisting of 
a principal $G$-bundle $E$ on $X_T$, together with a relative
canonical reduction $(L,f)$ of type $\tau$. Forgetting the reduction $(L,f)$
gives a $1$-morphism $i_{\tau} : Bun_{X/S}^{\tau}(G) \to Bun_{X/S}(G)$. 

It can be directly seen that the above morphism is representable and 
so $Bun_{X/S}^{\tau}(G)$ is an Artin stack,
even when the characteristic of $k$ is arbitrary and when the Behrend
conjecture does not hold (see [Gu-Ni 3]). In our present context, the 
Theorem \ref{main theorem}, together with 
Proposition \ref{semicontinuity and uniqueness}, 
show that the $1$-morphism $i_{\tau} : Bun_{X/S}^{\tau}(G) \to Bun_{X/S}(G)$
is a locally closed embedding of stacks 
(see [La-MB] Definition 3.14) (in particular, this 
gives another proof that $Bun_{X/S}^{\tau}(G)$)
is algebraic). This completes the proof of Theorem \ref{stacky main theorem}.

\bigskip 

\bigskip

\parskip=2pt

{\large \bf References}

\bigskip

[A-B] Atiyah, M.F. and Bott, R. : The Yang-Mills equations over 
Riemann surfaces. 
Philos. Trans. Roy. Soc. London A 308 (1983), 523-615.

[Be 1]  Behrend, K. : 
The Lefschetz Trace Formula for the Moduli Stack of
Principal Bundles. PhD Thesis, Berkeley. 
https://www.math.ubc.ca/~behrend/preprints.html

[Be 2] Behrend, K. : Semi-stability of reductive group schemes over curves. 
Math. Ann. 301 (1995), 281-305.

[EGA] Grothendieck, A. and Dieudonn\'e, J. 
: {\it  \'El\'ements de g\'eom\'etrie alg\'ebrique}, 
Publ. Math. IHES., vols. 4, 8, 11, 17, 20, 24, 28, 32 (1960-1967). 

[Gu-Ni 1] Gurjar, S. and Nitsure, N. : Schematic Harder-Narasimhan 
stratification for families of principal bundles and $\Lambda$-modules. 
Proc. Indian Acad. Sci. (Math. Sci.) 124 (2014), 315-332. 

[Gu-Ni 2] Gurjar, S. and Nitsure, N. : Schematic Harder-Narasimhan 
stratification for families of principal bundles in higher dimensions.
Preprint,  arXiv:1505.02236

[Gu-Ni 3] Gurjar, S. and Nitsure, N. : 
Harder-Narasimhan stacks for principal bundles in higher dimensions.
Preprint, arXiv:1605.08997.

[H-N] Harder, G. and Narasimhan, M. S. : 
On the cohomology groups of moduli spaces of vector bundles on curves. 
Math. Ann. 212 (1974/75), 215-248.

[H-AG] Hartshorne, R. : {\it Algebraic Geometry}, Springer (1977).

[La-MB] Laumon, G. and Moret-Bailly, L. : {\it Champs alg\'ebriques},
Springer (2000).

[K] Kleiman, S. : The Picard scheme. Part 5 
of {\it Fundamental Algebraic Geometry -- Grothendieck's FGA 
Explained}, Fantechi et al, Math. Surveys and Monographs Vol.
123, American Math. Soc. (2005).

[M-R] Mehta, V. B. and Ramanathan, A. : Semistable sheaves on 
projective varieties and their restriction to curves.
Math. Ann. 258 (1981/82), 213-224.

[Mu] Mumford, D. : {\it Lectures on Curves on an Algebraic Surface}, 
Annals of Mathematics Studies 59, Princeton University Press (1966)

[Ni 1] Nitsure, N. : Construction of Hilbert and Quot schemes. Part 2 
of {\it Fundamental Algebraic Geometry -- Grothendieck's FGA
Explained},  Fantechi et al, Math. Surveys and Monographs Vol.
123, American Math. Soc. (2005).

[Ni 2] Nitsure, N. : Deformation theory for vector bundles. Chapter 5
of {\it Moduli Spaces and Vector Bundles} 
(edited by Brambila-Paz, Bradlow, Garcia-Prada and Ramanan), 
London Math. Soc. Lect. Notes 359,  
Cambridge Univ. Press (2009).

[Ni 3] Nitsure, N. : 
Schematic Harder-Narasimhan stratification. 
Internat. J. Math. 22 (2011), 1365-1373. 

[Rag] Raghunathan, M. S. :
A note on quotients of real algebraic groups by arithmetic subgroups.
Invent. Math. 4 (1967/1968) 318-335.

[Ram] Ramanathan, A. : Moduli for principal bundles over algebraic curves.
Part I - Proc. Indian Acad. Sci. (Math. Sci.), 106(3) (1997), 301-328. 
Part II - Proc. Indian Acad. Sci. (Math. Sci.), 106(4) (1997), 421-449.

[Sh] Shatz, S.S. : The decomposition and specialization of
algebraic families of vector bundles. Compositio Math. 35
(1977), 163-187.

[Si] Simpson, C. Moduli of representations of the fundamental 
group of a smooth projective variety -I.   
Publ. Math. IHES 79 (1994), 47-129.

\bigskip

\bigskip

\bigskip

{\footnotesize

Sudarshan Gurjar \hfill Nitin Nitsure\\
Department of Mathematics  \hfill School of Mathematics\\
Indian Institute for Technology, Bombay 
\hfill Tata Institute of Fundamental Research\\
Powai \hfill Homi Bhabha Road\\
Mumbai 400 076 \hfill Mumbai 400 005\\ 
India \hfill India\\
{\tt srgurjar1984@gmail.com} \hfill {\tt nitsure@math.tifr.res.in}

\bigskip

\bigskip

}

\end{document}